\documentclass{amsart}
\usepackage{amsfonts}
\usepackage{graphicx}
\usepackage{amscd}

\newtheorem{theorem}{Theorem}
\theoremstyle{plain}

\newtheorem{corollary}{Corollary}

\newtheorem{lemma}{Lemma}

\newtheorem{remark}{Remark}

\numberwithin{equation}{section}

\input{tcilatex}

\begin{document}
\title[Gr\"{u}ss Type Inequalities for the Stieltjes Integral]{New
Inequalities of Gr\"{u}ss Type for the Stieltjes Integral and Applications}
\author{S.S. Dragomir}
\address{School of Communications and Informatics\\
Victoria University of Technology\\
P.O. Box 14428, MCMC 8001,\\
Victoria, Australia.}
\email{sever@matilda.vu.edu.au}
\urladdr{http://rgmia.vu.edu.au/SSDragomirWeb.html}
\date{June 18, 2002.}

\begin{abstract}
Sharp bounds of two \v{C}eby\v{s}ev functionals for the Stieltjes integrals
and applications for quadrature rules are given.
\end{abstract}

\maketitle

\section{Introduction}

Consider the \textit{weighted \v{C}eby\v{s}ev functional} 
\begin{multline}
T_{w}\left( f,g\right) :=\frac{1}{\int_{a}^{b}w\left( t\right) dt}%
\int_{a}^{b}w\left( t\right) f\left( t\right) g\left( t\right) dt
\label{1.1} \\
-\frac{1}{\int_{a}^{b}w\left( t\right) dt}\int_{a}^{b}w\left( t\right)
f\left( t\right) dt\cdot \frac{1}{\int_{a}^{b}w\left( t\right) dt}%
\int_{a}^{b}w\left( t\right) g\left( t\right) dt
\end{multline}
where $f,g,w:\left[ a,b\right] \rightarrow \mathbb{R}$ and $w\left( t\right)
\geq 0$ for a.e. $t\in \left[ a,b\right] $ are measurable functions such
that the involved integrals exist and $\int_{a}^{b}w\left( t\right) dt>0.$

In \cite{CD}, the authors obtained, among others, the following
inequalities: 
\begin{align}
& \left| T_{w}\left( f,g\right) \right|  \label{1.2} \\
& \leq \frac{1}{2}\left( M-m\right) \frac{1}{\int_{a}^{b}w\left( t\right) dt}%
\int_{a}^{b}w\left( t\right) \left| g\left( t\right) -\frac{1}{%
\int_{a}^{b}w\left( s\right) ds}\int_{a}^{b}w\left( s\right) g\left(
s\right) ds\right| dt  \notag \\
& \leq \frac{1}{2}\left( M-m\right) \left[ \frac{1}{\int_{a}^{b}w\left(
t\right) dt}\int_{a}^{b}w\left( t\right) \right.  \notag \\
& \;\;\;\;\;\;\;\;\;\;\;\;\;\;\;\;\;\;\;\;\;\times \left. \left| g\left(
t\right) -\frac{1}{\int_{a}^{b}w\left( s\right) ds}\int_{a}^{b}w\left(
s\right) g\left( s\right) ds\right| ^{p}dt\right] ^{\frac{1}{p}}\;\;\left(
p>1\right)  \notag \\
& \leq \frac{1}{2}\left( M-m\right) ess\sup\limits_{t\in \left[ a,b\right]
}\left| g\left( t\right) -\frac{1}{\int_{a}^{b}w\left( s\right) ds}%
\int_{a}^{b}w\left( s\right) g\left( s\right) ds\right|  \notag
\end{align}
provided 
\begin{equation}
-\infty <m\leq f\left( t\right) \leq M<\infty \text{ \hspace{0.05in}for a.e. 
}t\in \left[ a,b\right]  \label{1.3}
\end{equation}
and the corresponding integrals are finite. The constant $\frac{1}{2}$ is
sharp in all the inequalities in (\ref{1.2}) in the sense that it cannot be
replaced by a smaller constant.

In addition, if 
\begin{equation}
-\infty <n\leq g\left( t\right) \leq N<\infty \text{ \hspace{0.05in}for a.e. 
}t\in \left[ a,b\right] ,  \label{1.4}
\end{equation}
then the following refinement of the celebrated Gr\"{u}ss inequality is
obtained: 
\begin{align}
& \left| T_{w}\left( f,g\right) \right|  \label{1.5} \\
& \leq \frac{1}{2}\left( M-m\right) \frac{1}{\int_{a}^{b}w\left( t\right) dt}%
\int_{a}^{b}w\left( t\right) \left| g\left( t\right) -\frac{1}{%
\int_{a}^{b}w\left( s\right) ds}\int_{a}^{b}w\left( s\right) g\left(
s\right) ds\right| dt  \notag \\
& \leq \frac{1}{2}\left( M-m\right) \left[ \frac{1}{\int_{a}^{b}w\left(
t\right) dt}\int_{a}^{b}w\left( t\right) \right.  \notag \\
& \;\;\;\;\;\;\;\;\;\;\;\;\;\;\;\;\;\;\;\;\;\;\hspace{0.05in}\times \left.
\left| g\left( t\right) -\frac{1}{\int_{a}^{b}w\left( s\right) ds}%
\int_{a}^{b}w\left( s\right) g\left( s\right) ds\right| ^{2}dt\right] ^{%
\frac{1}{2}}  \notag \\
& \leq \frac{1}{4}\left( M-m\right) \left( N-n\right) .  \notag
\end{align}
Here, the constants $\frac{1}{2}$ and $\frac{1}{4}$ are also sharp in the
sense mentioned above.

In this paper, we extend the above results for Riemann-Stieltjes integrals.
A quadrature formula is also considered.

For this purpose, we introduce the following \v{C}eby\v{s}ev functional for
the Stieltjes integral 
\begin{multline}
T\left( f,g;u\right) :=\frac{1}{u\left( b\right) -u\left( a\right) }%
\int_{a}^{b}f\left( t\right) g\left( t\right) du\left( t\right)  \label{1.6}
\\
-\frac{1}{u\left( b\right) -u\left( a\right) }\int_{a}^{b}f\left( t\right)
du\left( t\right) \cdot \frac{1}{u\left( b\right) -u\left( a\right) }%
\int_{a}^{b}g\left( t\right) du\left( t\right) ,
\end{multline}
where $f,g\in C\left[ a,b\right] $ (are continuous on $\left[ a,b\right] $)
and $u\in BV\left[ a,b\right] $ (is of bounded variation on $\left[ a,b%
\right] $) with $u\left( b\right) \neq u\left( a\right) .$

For some recent inequalities for Stieltjes integral see \cite{D0}-\cite{D2}.

\section{Some Inequalities by Generalised \v{C}eby\v{s}ev Functional}

The following result holds \cite{6b}.

\begin{theorem}
\label{t2.1a}Let $f,g:\left[ a,b\right] \rightarrow \mathbb{R}$ be
continuous on $\left[ a,b\right] $ and $u:\left[ a,b\right] \rightarrow 
\mathbb{R}$ with $u\left( a\right) \neq u\left( b\right) .$ Assume also that
there exists the real constants $m,M$ such that 
\begin{equation}
m\leq f\left( t\right) \leq M\text{ \hspace{0.05in}for each }t\in \left[ a,b%
\right] .  \label{2.1a}
\end{equation}
If $u$ is of bounded variation on $\left[ a,b\right] ,$ then we have the
inequality 
\begin{multline}
\left| T\left( f,g;u\right) \right| \leq \frac{1}{2}\left( M-m\right) \frac{1%
}{\left| u\left( b\right) -u\left( a\right) \right| }  \label{2.2a} \\
\times \left\| g-\frac{1}{u\left( b\right) -u\left( a\right) }%
\int_{a}^{b}g\left( s\right) du\left( s\right) \right\| _{\infty
}\bigvee_{a}^{b}\left( u\right) ,
\end{multline}
where $\bigvee_{a}^{b}\left( u\right) $ denotes the total variation of $u$
in $\left[ a,b\right] .$ The constant $\frac{1}{2}$ is sharp, in the sense
that it cannot be replaced by a smaller constant.
\end{theorem}

\begin{proof}
It is easy to see, by simple computation with the Stieltjes integral, that
the following equality 
\begin{multline}
T\left( f,g;u\right) =\frac{1}{u\left( b\right) -u\left( a\right) }%
\int_{a}^{b}\left[ f\left( t\right) -\frac{m+M}{2}\right]  \label{2.3a} \\
\times \left[ g\left( t\right) -\frac{1}{u\left( b\right) -u\left( a\right) }%
\int_{a}^{b}g\left( s\right) du\left( s\right) \right] du\left( t\right)
\end{multline}
holds.

Using the known inequality 
\begin{equation}
\left| \int_{a}^{b}p\left( t\right) dv\left( t\right) \right| \leq
\sup\limits_{t\in \left[ a,b\right] }\left| p\left( t\right) \right|
\bigvee_{a}^{b}\left( v\right) ,  \label{2.4a}
\end{equation}
provided $p\in C\left[ a,b\right] $ and $v\in BV\left[ a,b\right] ,$ we
have, by (\ref{2.3a}), that 
\begin{align*}
\left| T\left( f,g;u\right) \right| & \leq \sup\limits_{t\in \left[ a,b%
\right] }\left| \left[ f\left( t\right) -\frac{m+M}{2}\right] \left[ g\left(
t\right) -\frac{1}{u\left( b\right) -u\left( a\right) }\int_{a}^{b}g\left(
s\right) du\left( s\right) \right] \right| \\
& \;\;\;\;\;\;\;\;\;\cdot \frac{1}{\left| u\left( b\right) -u\left( a\right)
\right| }\bigvee_{a}^{b}\left( u\right) \\
& \left( \text{since }\left| f\left( t\right) -\frac{m+M}{2}\right| \leq 
\frac{M-m}{2}\text{ for any }t\in \left[ a,b\right] \right) \\
& \leq \frac{M-m}{2}\left\| g-\frac{1}{u\left( b\right) -u\left( a\right) }%
\int_{a}^{b}g\left( s\right) du\left( s\right) \right\| _{\infty }\cdot 
\frac{1}{\left| u\left( b\right) -u\left( a\right) \right| }%
\bigvee_{a}^{b}\left( u\right)
\end{align*}
and the inequality (\ref{2.2a}) is proved.

To prove the sharpness of the constant $\frac{1}{2}$ in the inequality (\ref%
{2.2a}), we assume that it holds with a constant $C>0,$ i.e., 
\begin{multline}
\left| T\left( f,g;u\right) \right| \leq C\left( M-m\right) \frac{1}{\left|
u\left( b\right) -u\left( a\right) \right| }  \label{2.5a} \\
\times \left\| g-\frac{1}{u\left( b\right) -u\left( a\right) }%
\int_{a}^{b}g\left( s\right) du\left( s\right) \right\| _{\infty
}\bigvee_{a}^{b}\left( u\right) .
\end{multline}
Let us consider the functions $f=g,$ $f:\left[ a,b\right] \rightarrow 
\mathbb{R}$, $f\left( t\right) =t,$ $t\in \left[ a,b\right] $ and $u:\left[
a,b\right] \rightarrow \mathbb{R}$ given by 
\begin{equation}
u\left( t\right) =\left\{ 
\begin{array}{ll}
-1 & \text{if \hspace{0.05in}}t=a, \\ 
&  \\ 
0 & \text{if \hspace{0.05in}}t\in \left( a,b\right) , \\ 
&  \\ 
1 & \text{if \hspace{0.05in}}t=b.%
\end{array}
\right.  \label{2.6a}
\end{equation}
Then $f,g$ are continuous on $\left[ a,b\right] ,$ $u$ is of bounded
variation on $\left[ a,b\right] $ and 
\begin{equation*}
\frac{1}{u\left( b\right) -u\left( a\right) }\int_{a}^{b}f\left( t\right)
g\left( t\right) du\left( t\right) =\frac{b^{2}+a^{2}}{2},
\end{equation*}
\begin{equation*}
\frac{1}{u\left( b\right) -u\left( a\right) }\int_{a}^{b}f\left( t\right)
du\left( t\right) =\frac{b+a}{2},
\end{equation*}
\begin{equation*}
\left\| g-\frac{1}{u\left( b\right) -u\left( a\right) }\int_{a}^{b}g\left(
s\right) du\left( s\right) \right\| _{\infty }=\sup\limits_{t\in \left[ a,b%
\right] }\left| t-\frac{a+b}{2}\right| =\frac{b-a}{2}
\end{equation*}
and 
\begin{equation*}
\bigvee_{a}^{b}\left( u\right) =2,\;\;M=b,\;\;m=a.
\end{equation*}
Inserting these values in (\ref{2.5a}), we get 
\begin{equation*}
\left| \frac{a^{2}+b^{2}}{2}-\frac{\left( a+b\right) ^{2}}{4}\right| \leq
C\left( b-a\right) \cdot \frac{1}{2}\cdot \frac{\left( b-a\right) }{2}\cdot
2,
\end{equation*}
giving $C\geq \frac{1}{2}$, and the theorem is thus proved.
\end{proof}

The corresponding result for monotonic function $u$ is incorporated in the
following theorem \cite{6b}.

\begin{theorem}
\label{t2.2}Assume that $f$ and $g$ are as in Theorem \ref{t2.1a}. If $u:%
\left[ a,b\right] \rightarrow \mathbb{R}$ is monotonic nondecreasing on $%
\left[ a,b\right] ,$ then one has the inequality: 
\begin{multline}
\left| T\left( f,g;u\right) \right| \leq \frac{1}{2}\left( M-m\right) \frac{1%
}{u\left( b\right) -u\left( a\right) }  \label{2.7a} \\
\times \int_{a}^{b}\left| g\left( t\right) -\frac{1}{u\left( b\right)
-u\left( a\right) }\int_{a}^{b}g\left( s\right) du\left( s\right) \right|
du\left( t\right) .
\end{multline}
The constant $\frac{1}{2}$ is sharp in the sense that it cannot be replaced
by a smaller constant.
\end{theorem}

\begin{proof}
Using the known inequality 
\begin{equation}
\left| \int_{a}^{b}p\left( t\right) dv\left( t\right) \right| \leq
\int_{a}^{b}\left| p\left( t\right) \right| dv\left( t\right) ,  \label{2.8a}
\end{equation}
provided $p\in C\left[ a,b\right] $ and $v$ is a monotonic nondecreasing
function on $\left[ a,b\right] ,$ we have (by the use of equality (\ref{2.3a}%
)) that 
\begin{align*}
\left| T\left( f,g;u\right) \right| & \leq \frac{1}{u\left( b\right)
-u\left( a\right) }\int_{a}^{b}\left| f\left( t\right) -\frac{m+M}{2}\right|
\\
& \;\;\;\;\;\;\;\;\;\;\;\;\;\;\;\;\;\;\;\;\;\;\times \left| g\left( t\right)
-\frac{1}{u\left( b\right) -u\left( a\right) }\int_{a}^{b}g\left( s\right)
du\left( s\right) \right| du\left( t\right) \\
& \leq \frac{1}{2}\left( M-m\right) \frac{1}{u\left( b\right) -u\left(
a\right) }\int_{a}^{b}\left| g\left( t\right) -\frac{1}{u\left( b\right)
-u\left( a\right) }\int_{a}^{b}g\left( s\right) du\left( s\right) \right|
du\left( t\right) .
\end{align*}
Now, assume that the inequality (\ref{2.7a}) holds with a constant $D>0,$
instead of $\frac{1}{2},$ i.e., 
\begin{multline}
\left| T\left( f,g;u\right) \right| \leq D\left( M-m\right) \frac{1}{u\left(
b\right) -u\left( a\right) }  \label{2.9ab} \\
\times \int_{a}^{b}\left| g\left( t\right) -\frac{1}{u\left( b\right)
-u\left( a\right) }\int_{a}^{b}g\left( s\right) du\left( s\right) \right|
du\left( t\right) .
\end{multline}
If we choose the same function as in the proof of Theorem \ref{t2.1a}, we
observe that $f,g$ are continuous and $u$ is monotonic nondecreasing on $%
\left[ a,b\right] .$ Then, for these functions, we have 
\begin{equation*}
T\left( f,g;u\right) =\frac{a^{2}+b^{2}}{2}-\frac{\left( a+b\right) ^{2}}{4}=%
\frac{\left( b-a\right) ^{2}}{4},
\end{equation*}
\begin{align*}
\int_{a}^{b}\left| g\left( t\right) -\frac{1}{u\left( b\right) -u\left(
a\right) }\int_{a}^{b}g\left( s\right) du\left( s\right) \right| du\left(
t\right) & =\int_{a}^{b}\left| t-\frac{a+b}{2}\right| du\left( t\right) \\
& =b-a,
\end{align*}
and then, by (\ref{2.9ab}) we get 
\begin{equation*}
\frac{\left( b-a\right) ^{2}}{4}\leq D\left( b-a\right) \frac{1}{2}\left(
b-a\right)
\end{equation*}
giving $D\geq \frac{1}{2}$, and the theorem is completely proved.
\end{proof}

The case when $u$ is a Lipschitzian function is embodied in the following
theorem \cite{6b}.

\begin{theorem}
\label{t2.3a}Assume that $f,g:\left[ a,b\right] \rightarrow \mathbb{R}$ are
Riemann integrable functions on $\left[ a,b\right] $ and $f$ satisfies the
condition (\ref{2.1a}). If $u:\left( a,b\right) \rightarrow \mathbb{R}$ $%
\left( u\left( b\right) \neq u\left( a\right) \right) $ is Lipschitzian with
the constant $L,$ then we have the inequality 
\begin{multline}
\left| T\left( f,g;u\right) \right| \leq \frac{1}{2}L\left( M-m\right) \frac{%
1}{\left| u\left( b\right) -u\left( a\right) \right| }  \label{2.10a} \\
\times \int_{a}^{b}\left| g\left( t\right) -\frac{1}{u\left( b\right)
-u\left( a\right) }\int_{a}^{b}g\left( s\right) du\left( s\right) \right| dt.
\end{multline}
The constant $\frac{1}{2}$ cannot be replaced by a smaller constant.
\end{theorem}

\begin{proof}
It is well known that if $p:\left[ a,b\right] \rightarrow \mathbb{R}$ is
Riemann integrable on $\left[ a,b\right] $ and $v:\left[ a,b\right]
\rightarrow \mathbb{R}$ is Lipschitzian with the constant $L,$ then the
Riemann-Stieltjes integral $\int_{a}^{b}p\left( t\right) dv\left( t\right) $
exists and 
\begin{equation}
\left| \int_{a}^{b}p\left( t\right) dv\left( t\right) \right| \leq
L\int_{a}^{b}\left| p\left( t\right) \right| dt.  \label{2.11a}
\end{equation}
Using this fact and the identity (\ref{2.3a}), we deduce 
\begin{align*}
\left| T\left( f,g;u\right) \right| & \leq \frac{L}{\left| u\left( b\right)
-u\left( a\right) \right| }\int_{a}^{b}\left| f\left( t\right) -\frac{m+M}{2}%
\right| \\
& \;\;\;\;\;\;\;\;\;\;\;\;\;\;\;\;\times \left| g\left( t\right) -\frac{1}{%
u\left( b\right) -u\left( a\right) }\int_{a}^{b}g\left( s\right) du\left(
s\right) \right| dt \\
& \leq \frac{1}{2}\left( M-m\right) \frac{L}{\left| u\left( b\right)
-u\left( a\right) \right| }\int_{a}^{b}\left| g\left( t\right) -\frac{1}{%
u\left( b\right) -u\left( a\right) }\int_{a}^{b}g\left( s\right) du\left(
s\right) \right| dt
\end{align*}
and the inequality (\ref{2.10a}) is proved.

Now, assume that (\ref{2.10a}) holds with a constant $E>0$ instead of $\frac{%
1}{2},$ i.e., 
\begin{multline}
\left| T\left( f,g;u\right) \right| \leq EL\left( M-m\right) \frac{1}{\left|
u\left( b\right) -u\left( a\right) \right| }  \label{2.12a} \\
\times \int_{a}^{b}\left| g\left( t\right) -\frac{1}{u\left( b\right)
-u\left( a\right) }\int_{a}^{b}g\left( s\right) du\left( s\right) \right| dt.
\end{multline}
Consider the function $f=g,$ $f:\left[ a,b\right] \rightarrow \mathbb{R}$
with 
\begin{equation*}
f\left( t\right) =\left\{ 
\begin{array}{ll}
-1 & \text{if \hspace{0.05in}}t\in \left[ a,\frac{a+b}{2}\right] \\ 
&  \\ 
1 & \text{if \hspace{0.05in}}t\in \left( \frac{a+b}{2},b\right]%
\end{array}
\right.
\end{equation*}
and $u:\left[ a,b\right] \rightarrow \mathbb{R}$, $u\left( t\right) =t.$
Then, obviously, $f$ and $g$ are Riemann integrable on $\left[ a,b\right] $
and $u$ is Lipschitzian with the constant $L=1.$

Since 
\begin{gather*}
\frac{1}{u\left( b\right) -u\left( a\right) }\int_{a}^{b}f\left( t\right)
g\left( t\right) du\left( t\right) =\frac{1}{b-a}\int_{a}^{b}dt=1, \\
\frac{1}{u\left( b\right) -u\left( a\right) }\int_{a}^{b}f\left( t\right)
du\left( t\right) =\frac{1}{u\left( b\right) -u\left( a\right) }%
\int_{a}^{b}g\left( t\right) du\left( t\right) =0, \\
\int_{a}^{b}\left\vert g\left( t\right) -\frac{1}{u\left( b\right) -u\left(
a\right) }\int_{a}^{b}g\left( s\right) du\left( s\right) \right\vert
dt=\int_{a}^{b}dt=b-a
\end{gather*}
and 
\begin{equation*}
M=1,\;\;m=1
\end{equation*}
then, by (\ref{2.12a}), we deduce $E\geq \frac{1}{2},$ and the theorem is
completely proved.
\end{proof}

The following result holds \cite{7b}.

\begin{theorem}
\label{t2.1}Let $f,g:\left[ a,b\right] \rightarrow \mathbb{R}$ be such that $%
f$ is of $r-H-$H\"{o}lder type on $\left[ a,b\right] ,$ i.e., 
\begin{equation}
\left| f\left( t\right) -f\left( s\right) \right| \leq H\left| t-s\right|
^{r}\text{ \hspace{0.05in}for any \hspace{0.05in}}t,s\in \left[ a,b\right] ,
\label{2.1}
\end{equation}
and $g$ is continuous on $\left[ a,b\right] .$ If $u:\left[ a,b\right]
\rightarrow \mathbb{R}$ is of bounded variation on $\left[ a,b\right] $ with 
$u\left( a\right) \neq u\left( b\right) ,$ then we have the inequality 
\begin{multline}
\left| T\left( f,g;u\right) \right| \leq \frac{H\left( b-a\right) ^{r}}{2^{r}%
}\cdot \frac{1}{\left| u\left( b\right) -u\left( a\right) \right| }
\label{2.2} \\
\times \left\| g-\frac{1}{u\left( b\right) -u\left( a\right) }%
\int_{a}^{b}g\left( s\right) du\left( s\right) \right\| _{\infty
}\bigvee_{a}^{b}\left( u\right) ,
\end{multline}
where $\bigvee_{a}^{b}\left( u\right) $ denotes the total variation of $u$
on $\left[ a,b\right] .$
\end{theorem}

\begin{proof}
It is easy to see, by simple computation with the Stieltjes integral, that
the following equality 
\begin{multline}
T\left( f,g;u\right) =\frac{1}{u\left( b\right) -u\left( a\right) }%
\int_{a}^{b}\left[ f\left( t\right) -f\left( \frac{a+b}{2}\right) \right]
\label{2.3} \\
\times \left[ g\left( t\right) -\frac{1}{u\left( b\right) -u\left( a\right) }%
\int_{a}^{b}g\left( s\right) du\left( s\right) \right] du\left( t\right)
\end{multline}
holds.

Using the known inequality 
\begin{equation}
\left| \int_{a}^{b}p\left( t\right) dv\left( t\right) \right| \leq
\sup\limits_{t\in \left[ a,b\right] }\left| p\left( t\right) \right|
\bigvee_{a}^{b}\left( v\right)  \label{2.4}
\end{equation}
provided $p\in C\left[ a,b\right] $ and $v\in BV\left[ a,b\right] ,$ we
have, by (\ref{2.3}), that 
\begin{align*}
\left| T\left( f,g;u\right) \right| & \leq \sup\limits_{t\in \left[ a,b%
\right] }\left| \left[ f\left( t\right) -f\left( \frac{a+b}{2}\right) \right]
\left[ g\left( t\right) -\frac{1}{u\left( b\right) -u\left( a\right) }%
\int_{a}^{b}g\left( s\right) du\left( s\right) \right] \right| \\
& \;\;\;\;\;\;\;\;\;\;\;\;\;\times \frac{1}{\left| u\left( b\right) -u\left(
a\right) \right| }\bigvee_{a}^{b}\left( u\right) \\
& \leq \sup\limits_{t\in \left[ a,b\right] }\left| f\left( t\right) -f\left( 
\frac{a+b}{2}\right) \right| \left\| g-\frac{1}{u\left( b\right) -u\left(
a\right) }\int_{a}^{b}g\left( s\right) du\left( s\right) \right\| _{\infty }
\\
& \;\;\;\;\;\;\;\;\;\;\;\;\;\times \frac{1}{\left| u\left( b\right) -u\left(
a\right) \right| }\bigvee_{a}^{b}\left( u\right) \\
& \leq L\left( \frac{b-a}{2}\right) ^{r}\left\| g-\frac{1}{u\left( b\right)
-u\left( a\right) }\int_{a}^{b}g\left( s\right) du\left( s\right) \right\|
_{\infty } \\
& \;\;\;\;\;\;\;\;\;\;\;\;\;\times \frac{1}{\left| u\left( b\right) -u\left(
a\right) \right| }\bigvee_{a}^{b}\left( u\right) ,
\end{align*}
and the inequality (\ref{2.2}) is proved.
\end{proof}

The following corollary may be useful in applications \cite{7b}.

\begin{corollary}
\label{c2.2}Let $f$ be Lipschitzian with the constant $L>0,$ i.e., 
\begin{equation}
\left\vert f\left( t\right) -f\left( s\right) \right\vert \leq L\left\vert
t-s\right\vert \text{ \hspace{0.05in}for any \hspace{0.05in}}t,s\in \left[
a,b\right] ,  \label{2.5}
\end{equation}
and $u,g$ are as in Theorem \ref{t2.1}. Then we have the inequality 
\begin{multline}
\left\vert T\left( f,g;u\right) \right\vert \leq \frac{1}{2}\frac{L\left(
b-a\right) }{\left\vert u\left( b\right) -u\left( a\right) \right\vert }
\label{2.6} \\
\times \left\Vert g-\frac{1}{u\left( b\right) -u\left( a\right) }%
\int_{a}^{b}g\left( s\right) du\left( s\right) \right\Vert _{\infty
}\bigvee_{a}^{b}\left( u\right) .
\end{multline}
The constant $\frac{1}{2}$ cannot be replaced by a smaller constant.
\end{corollary}

\begin{proof}
The inequality (\ref{2.6}) follows by (\ref{2.2}) for $r=1.$ It remains to
prove only the sharpness of the constant $\frac{1}{2}.$

Consider the functions $f=g,$ where $f:\left[ a,b\right] \rightarrow \mathbb{%
R}$, $f\left( t\right) =t$ and $u:\left[ a,b\right] \rightarrow \mathbb{R}$,
given by 
\begin{equation}
u\left( t\right) =\left\{ 
\begin{array}{ll}
-1 & \text{if \hspace{0.05in}}t=a, \\ 
&  \\ 
0 & \text{if \hspace{0.05in}}t\in \left( a,b\right) , \\ 
&  \\ 
1 & \text{if \hspace{0.05in}}t=b.%
\end{array}
\right.  \label{2.6.a}
\end{equation}
Then, $f$ is Lipschitzian with the constant $L=1,$ $g$ is continuous and $u$
is of bounded variation.

If we assume that the inequality (\ref{2.6}) holds with a constant $C>0,$
i.e., 
\begin{equation}
\left| T\left( f,g;u\right) \right| \leq CL\left( b-a\right) \left\| g-\frac{%
1}{u\left( b\right) -u\left( a\right) }\int_{a}^{b}g\left( s\right) du\left(
s\right) \right\| _{\infty }\bigvee_{a}^{b}\left( u\right) ,  \label{2.7}
\end{equation}
and since 
\begin{equation*}
\frac{1}{u\left( b\right) -u\left( a\right) }\int_{a}^{b}f\left( t\right)
g\left( t\right) du\left( t\right) =\frac{b^{2}+a^{2}}{2},
\end{equation*}
\begin{equation*}
\frac{1}{u\left( b\right) -u\left( a\right) }\int_{a}^{b}f\left( t\right)
du\left( t\right) =\frac{1}{u\left( b\right) -u\left( a\right) }%
\int_{a}^{b}g\left( t\right) du\left( t\right) =\frac{b+a,}{2}
\end{equation*}
\begin{equation*}
\left\| g-\frac{1}{u\left( b\right) -u\left( a\right) }\int_{a}^{b}g\left(
s\right) du\left( s\right) \right\| _{\infty }=\sup\limits_{t\in \left[ a,b%
\right] }\left| t-\frac{a+b}{2}\right| =\frac{b-a}{2}
\end{equation*}
and $\bigvee_{a}^{b}\left( u\right) =2,$ then, by (\ref{2.7}), we have 
\begin{equation*}
\left| \frac{b^{2}+a^{2}}{2}-\left( \frac{a+b}{2}\right) ^{2}\right| \leq C%
\frac{\left( b-a\right) }{2}\frac{b-a}{2}\cdot 2,
\end{equation*}
giving $C\geq \frac{1}{2}.$
\end{proof}

The following result concerning monotonic function $u:\left[ a,b\right]
\rightarrow \mathbb{R}$ also holds \cite{7b}.

\begin{theorem}
\label{t2.3}Assume that $f$ and $g$ are as in Theorem \ref{t2.1}. If $u:%
\left[ a,b\right] \rightarrow \mathbb{R}$ is monotonic nondecreasing on $%
\left[ a,b\right] $ with $u\left( b\right) >u\left( a\right) ,$ then we have
the inequalities: 
\begin{align}
\left| T\left( f,g;u\right) \right| & \leq \frac{H}{u\left( b\right)
-u\left( a\right) }\int_{a}^{b}\left| t-\frac{a+b}{2}\right| ^{r}
\label{2.8} \\
& \;\;\;\;\;\;\;\times \left| g\left( t\right) -\frac{1}{u\left( b\right)
-u\left( a\right) }\int_{a}^{b}g\left( s\right) du\left( s\right) \right|
du\left( t\right)  \notag \\
& \leq \frac{H\left( b-a\right) ^{r}}{2^{r}\left[ u\left( b\right) -u\left(
a\right) \right] }  \notag \\
& \;\;\;\;\;\;\;\;\times \int_{a}^{b}\left| g\left( t\right) -\frac{1}{%
u\left( b\right) -u\left( a\right) }\int_{a}^{b}g\left( s\right) du\left(
s\right) \right| du\left( t\right) .  \notag
\end{align}
\end{theorem}

\begin{proof}
Using the known inequality 
\begin{equation}
\left| \int_{a}^{b}p\left( t\right) dv\left( t\right) \right| \leq
\int_{a}^{b}\left| p\left( t\right) \right| dv\left( t\right) ,  \label{2.9}
\end{equation}
provided $p\in C\left[ a,b\right] $ and $v$ is monotonic nondecreasing on $%
\left[ a,b\right] ,$ we have, by (\ref{2.3}), the following estimate: 
\begin{align*}
\left| T\left( f,g;u\right) \right| & \leq \frac{1}{u\left( b\right)
-u\left( a\right) }\int_{a}^{b}\left| \left( f\left( t\right) -f\left( \frac{%
a+b}{2}\right) \right) \right. \\
& \;\;\;\times \left. \left( g\left( t\right) -\frac{1}{u\left( b\right)
-u\left( a\right) }\int_{a}^{b}g\left( s\right) du\left( s\right) \right)
\right| du\left( t\right) \\
& \leq \frac{H}{u\left( b\right) -u\left( a\right) }\int_{a}^{b}\left| t-%
\frac{a+b}{2}\right| ^{r} \\
& \;\;\;\times \left| g\left( t\right) -\frac{1}{u\left( b\right) -u\left(
a\right) }\int_{a}^{b}g\left( s\right) du\left( s\right) \right| du\left(
t\right) \\
& \leq \frac{H}{u\left( b\right) -u\left( a\right) }\sup\limits_{t\in \left[
a,b\right] }\left| t-\frac{a+b}{2}\right| ^{r} \\
& \;\;\;\;\;\times \int_{a}^{b}\left| g\left( t\right) -\frac{1}{u\left(
b\right) -u\left( a\right) }\int_{a}^{b}g\left( s\right) du\left( s\right)
\right| du\left( t\right)
\end{align*}
which simply provides (\ref{2.8}).
\end{proof}

The particular case of Lipschitzian functions that is relevant for
applications is embodied in the following corollary \cite{7b}.

\begin{corollary}
\label{c2.4}Assume that $f$ is $L-$Lipschitzian, $g$ is continuous and $u$
is monotonic nondecreasing on $\left[ a,b\right] $ with $u\left( b\right)
>u\left( a\right) .$ Then we have the inequalities 
\begin{eqnarray}
\left| T\left( f,g;u\right) \right| &\leq &\frac{L}{u\left( b\right)
-u\left( a\right) }\int_{a}^{b}\left| t-\frac{a+b}{2}\right|  \label{2.9a} \\
&&\;\;\;\times \left| g\left( t\right) -\frac{1}{u\left( b\right) -u\left(
a\right) }\int_{a}^{b}g\left( s\right) du\left( s\right) \right| du\left(
t\right)  \notag \\
&\leq &\frac{1}{2}\cdot \frac{L\left( b-a\right) }{u\left( b\right) -u\left(
a\right) }  \notag \\
&&\;\;\;\times \int_{a}^{b}\left| g\left( t\right) -\frac{1}{u\left(
b\right) -u\left( a\right) }\int_{a}^{b}g\left( s\right) du\left( s\right)
\right| du\left( t\right) .  \notag
\end{eqnarray}
The first inequality is sharp. The constant $\frac{1}{2}$ in the second
inequality cannot be replaced by a smaller constant.
\end{corollary}

\begin{proof}
The inequality (\ref{2.9a}) follows by (\ref{2.8}) on choosing $r=1.$ Assume
that (\ref{2.9a}) holds with the constants $D,E>0,$ i.e., 
\begin{align}
& \left| T\left( f,g;u\right) \right|  \label{2.10} \\
& \leq \frac{LD}{u\left( b\right) -u\left( a\right) }\int_{a}^{b}\left| t-%
\frac{a+b}{2}\right|  \notag \\
& \;\;\;\;\;\;\;\times \left| g\left( t\right) -\frac{1}{u\left( b\right)
-u\left( a\right) }\int_{a}^{b}g\left( s\right) du\left( s\right) \right|
du\left( t\right)  \notag \\
& \leq \frac{LE\left( b-a\right) }{u\left( b\right) -u\left( a\right) }%
\int_{a}^{b}\left| g\left( t\right) -\frac{1}{u\left( b\right) -u\left(
a\right) }\int_{a}^{b}g\left( s\right) du\left( s\right) \right| du\left(
t\right) .  \notag
\end{align}
Consider the functions $f=g$, where $f:\left[ a,b\right] \rightarrow \mathbb{%
R}$, $f\left( t\right) =t$ and $u$ is as given by (\ref{2.6.a}). Then,
obviously, $f$ is Lipschitzian with the constant $L=1,$ $g$ is continuous
and $u$ is monotonic nondecreasing on $\left[ a,b\right] .$

Since, we know, for these functions 
\begin{equation*}
T\left( f,g;u\right) =\frac{\left( b-a\right) ^{2}}{4},
\end{equation*}
and 
\begin{equation*}
\int_{a}^{b}\left| t-\frac{a+b}{2}\right| \left| g\left( t\right) -\frac{1}{%
u\left( b\right) -u\left( a\right) }\int_{a}^{b}g\left( s\right) du\left(
s\right) \right| du\left( t\right) =\frac{\left( b-a\right) ^{2}}{2},
\end{equation*}
\begin{equation*}
\int_{a}^{b}\left| g\left( t\right) -\frac{1}{u\left( b\right) -u\left(
a\right) }\int_{a}^{b}g\left( s\right) du\left( s\right) \right| du\left(
t\right) =b-a,
\end{equation*}
then by (\ref{2.10}) we deduce 
\begin{equation*}
\frac{\left( b-a\right) ^{2}}{4}\leq \frac{D}{2}\cdot \frac{\left(
b-a\right) ^{2}}{2}\leq \frac{E\left( b-a\right) ^{2}}{2}
\end{equation*}
giving $D\geq 1$ and $E\geq \frac{1}{2}.$
\end{proof}

Another natural possibility to obtain bounds for the functional $T\left(
f,g;u\right) ,$ where $u$ is Lipschitzian with the constant $K>0,$ is
embodied in the following theorem \cite{7b}.

\begin{theorem}
\label{t2.5}Assume that $f:\left[ a,b\right] \rightarrow \mathbb{R}$ is of $%
r-H-$H\"{o}lder type on $\left[ a,b\right] .$ If $g:\left[ a,b\right]
\rightarrow \mathbb{R}$ is Riemann integrable on $\left[ a,b\right] $ and $u:%
\left[ a,b\right] \rightarrow \mathbb{R}$ is Lipschitzian with the constant $%
K>0$ and $u\left( a\right) \neq u\left( b\right) ,$ then one has the
inequalities: 
\begin{align}
& \left| T\left( f,g;u\right) \right|  \label{2.11} \\
& \leq \frac{HK}{\left| u\left( b\right) -u\left( a\right) \right| }%
\int_{a}^{b}\left| t-\frac{a+b}{2}\right| ^{r}  \notag \\
& \;\;\;\;\;\times \left| g\left( t\right) -\frac{1}{u\left( b\right)
-u\left( a\right) }\int_{a}^{b}g\left( s\right) du\left( s\right) \right| dt
\notag \\
& \leq \left\{ 
\begin{array}{l}
\frac{HK\left( b-a\right) ^{r+1}}{2^{r}\left( r+1\right) \left| u\left(
b\right) -u\left( a\right) \right| }\left\| g-\frac{1}{u\left( b\right)
-u\left( a\right) }\int_{a}^{b}g\left( s\right) du\left( s\right) \right\|
_{\infty }; \\ 
\\ 
\frac{HK\left( b-a\right) ^{r+\frac{1}{q}}}{2^{r}\left( qr+1\right) ^{\frac{1%
}{q}}\left| u\left( b\right) -u\left( a\right) \right| }\left\| g-\frac{1}{%
u\left( b\right) -u\left( a\right) }\int_{a}^{b}g\left( s\right) du\left(
s\right) \right\| _{p}\text{ \hspace{0.05in}} \\ 
\hfill \text{if \hspace{0.05in}}p>1,\;\frac{1}{p}+\frac{1}{q}=1; \\ 
\\ 
\frac{HK\left( b-a\right) ^{r}}{2^{r}\left| u\left( b\right) -u\left(
a\right) \right| }\left\| g-\frac{1}{u\left( b\right) -u\left( a\right) }%
\int_{a}^{b}g\left( s\right) du\left( s\right) \right\| _{1}.%
\end{array}
\right.  \notag
\end{align}
\end{theorem}

\begin{proof}
Using the identity (\ref{2.3}), we have successively 
\begin{align}
\left| T\left( f,g;u\right) \right| & \leq \frac{K}{\left| u\left( b\right)
-u\left( a\right) \right| }\int_{a}^{b}\left| f\left( t\right) -f\left( 
\frac{a+b}{2}\right) \right|  \label{2.12} \\
& \;\;\;\;\;\;\;\;\;\;\;\;\;\;\times \left| g\left( t\right) -\frac{1}{%
u\left( b\right) -u\left( a\right) }\int_{a}^{b}g\left( s\right) du\left(
s\right) \right| dt  \notag \\
& \leq \frac{KH}{\left| u\left( b\right) -u\left( a\right) \right| }%
\int_{a}^{b}\left| t-\frac{a+b}{2}\right| ^{r}  \notag \\
& \;\;\;\;\;\;\;\;\;\;\;\;\;\;\;\times \left| g\left( t\right) -\frac{1}{%
u\left( b\right) -u\left( a\right) }\int_{a}^{b}g\left( s\right) du\left(
s\right) \right| dt  \notag
\end{align}
and the first inequality in (\ref{2.11}) is proved.

Since 
\begin{align*}
& \int_{a}^{b}\left| t-\frac{a+b}{2}\right| ^{r}\left| g\left( t\right) -%
\frac{1}{u\left( b\right) -u\left( a\right) }\int_{a}^{b}g\left( s\right)
du\left( s\right) \right| dt \\
& \leq \left\| g-\frac{1}{u\left( b\right) -u\left( a\right) }%
\int_{a}^{b}g\left( s\right) du\left( s\right) \right\| _{\infty
}\int_{a}^{b}\left| t-\frac{a+b}{2}\right| ^{r}dt \\
& =\frac{\left( b-a\right) ^{r+1}}{2^{r}\left( r+1\right) }\left\| g-\frac{1%
}{u\left( b\right) -u\left( a\right) }\int_{a}^{b}g\left( s\right) du\left(
s\right) \right\| _{\infty },
\end{align*}
then by (\ref{2.12}) we deduce the first part in the second inequality in (%
\ref{2.11}).

By H\"{o}lder's integral inequality we have 
\begin{align*}
& \int_{a}^{b}\left| t-\frac{a+b}{2}\right| ^{r}\left| g\left( t\right) -%
\frac{1}{u\left( b\right) -u\left( a\right) }\int_{a}^{b}g\left( s\right)
du\left( s\right) \right| dt \\
& \leq \left( \int_{a}^{b}\left| t-\frac{a+b}{2}\right| ^{qr}dt\right) ^{%
\frac{1}{q}}\left( \int_{a}^{b}\left| g\left( t\right) -\frac{1}{u\left(
b\right) -u\left( a\right) }\int_{a}^{b}g\left( s\right) du\left( s\right)
\right| ^{p}dt\right) ^{\frac{1}{p}} \\
& =\left[ \frac{\left( b-a\right) ^{qr+1}}{2^{qr}\left( qr+1\right) }\right]
^{\frac{1}{q}}\left\| g-\frac{1}{u\left( b\right) -u\left( a\right) }%
\int_{a}^{b}g\left( s\right) du\left( s\right) \right\| _{p} \\
& =\frac{\left( b-a\right) ^{r+\frac{1}{q}}}{2^{r}\left( qr+1\right) ^{\frac{%
1}{q}}}\left\| g-\frac{1}{u\left( b\right) -u\left( a\right) }%
\int_{a}^{b}g\left( s\right) du\left( s\right) \right\| _{p}.
\end{align*}
Using (\ref{2.12}), we deduce the second part of the second inequality in (%
\ref{2.11}).

Finally, since 
\begin{equation*}
\left| t-\frac{a+b}{2}\right| ^{r}\leq \left( \frac{b-a}{2}\right)
^{r},\;\;\;t\in \left[ a,b\right] ,
\end{equation*}
we deduce 
\begin{multline*}
\int_{a}^{b}\left| t-\frac{a+b}{2}\right| ^{r}\left| g\left( t\right) -\frac{%
1}{u\left( b\right) -u\left( a\right) }\int_{a}^{b}g\left( s\right) du\left(
s\right) \right| dt \\
\leq \frac{\left( b-a\right) ^{r}}{2^{r}}\left\| g-\frac{1}{u\left( b\right)
-u\left( a\right) }\int_{a}^{b}g\left( s\right) du\left( s\right) \right\|
_{1}
\end{multline*}
and the theorem is completely proved.
\end{proof}

The following particular case is useful in applications \cite{7b}.

\begin{corollary}
\label{c2.6}If $f$ is Lipschitzian with the constant $L$ and $g$ and $u$ are
as in Theorem \ref{t2.5}, then we have the inequalities: 
\begin{multline}
\left| T\left( f,g;u\right) \right| \leq \frac{LK}{\left| u\left( b\right)
-u\left( a\right) \right| }\int_{a}^{b}\left| t-\frac{a+b}{2}\right|
\label{2.13} \\
\times \left| g\left( t\right) -\frac{1}{u\left( b\right) -u\left( a\right) }%
\int_{a}^{b}g\left( s\right) du\left( s\right) \right| dt
\end{multline}
\begin{equation*}
\leq \left\{ 
\begin{array}{l}
\dfrac{LK\left( b-a\right) ^{2}}{4\left| u\left( b\right) -u\left( a\right)
\right| }\left\| g-\dfrac{1}{u\left( b\right) -u\left( a\right) }%
\int_{a}^{b}g\left( s\right) du\left( s\right) \right\| _{\infty }; \\ 
\\ 
\dfrac{LK\left( b-a\right) ^{1+\frac{1}{q}}}{2\left( q+1\right) ^{\frac{1}{q}%
}\left| u\left( b\right) -u\left( a\right) \right| }\left\| g-\dfrac{1}{%
u\left( b\right) -u\left( a\right) }\int_{a}^{b}g\left( s\right) du\left(
s\right) \right\| _{p}\text{ \hspace{0.05in}} \\ 
\hfill \text{if \hspace{0.05in}}p>1,\;\frac{1}{p}+\frac{1}{q}=1; \\ 
\\ 
\dfrac{LK\left( b-a\right) }{2\left| u\left( b\right) -u\left( a\right)
\right| }\left\| g-\dfrac{1}{u\left( b\right) -u\left( a\right) }%
\int_{a}^{b}g\left( s\right) du\left( s\right) \right\| _{1}.%
\end{array}
\right.
\end{equation*}
The first inequality in (\ref{2.13}) is sharp.

The constants $\frac{1}{4}$ and $\frac{1}{2}$ in the second branch of the
second inequality cannot be replaced by smaller constants, respectively.
\end{corollary}

\begin{proof}
The inequality (\ref{2.13}) follows obviously from (\ref{2.11}) on choosing $%
r=1.$

Now, assume that the following inequalities hold 
\begin{align}
& \left| T\left( f,g;u\right) \right|  \label{2.14} \\
& \leq \frac{CLK}{\left| u\left( b\right) -u\left( a\right) \right| }%
\int_{a}^{b}\left| t-\frac{a+b}{2}\right|  \notag \\
& \;\;\;\;\;\times \left| g\left( t\right) -\frac{1}{u\left( b\right)
-u\left( a\right) }\int_{a}^{b}g\left( s\right) du\left( s\right) \right| dt
\notag \\
& \leq \left\{ 
\begin{array}{l}
\dfrac{DLK\left( b-a\right) ^{2}}{\left| u\left( b\right) -u\left( a\right)
\right| }\left\| g-\dfrac{1}{u\left( b\right) -u\left( a\right) }%
\dint_{a}^{b}g\left( s\right) du\left( s\right) \right\| _{\infty }; \\ 
\\ 
\dfrac{ELK\left( b-a\right) ^{1+\frac{1}{q}}}{\left( q+1\right) ^{\frac{1}{q}%
}\left| u\left( b\right) -u\left( a\right) \right| }\left\| g-\dfrac{1}{%
u\left( b\right) -u\left( a\right) }\dint_{a}^{b}g\left( s\right) du\left(
s\right) \right\| _{p} \\ 
\hfill \text{if \hspace{0.05in}}p>1,\;\frac{1}{p}+\frac{1}{q}=1;%
\end{array}
\right.  \notag
\end{align}
with $C,D,E>0.$

Consider the functions $f,g,u:\left[ a,b\right] \rightarrow \mathbb{R}$,
defined by $f\left( t\right) =t-\frac{a+b}{2},$ $u\left( t\right) =t$ and 
\begin{equation*}
g\left( t\right) =\left\{ 
\begin{array}{ll}
-1 & \text{if \hspace{0.05in}}t\in \left[ a,\frac{a+b}{2}\right] , \\ 
&  \\ 
1 & \text{if \hspace{0.05in}}t\in \left( \frac{a+b}{2},b\right] .%
\end{array}
\right.
\end{equation*}
Then both $f$ and $u$ are Lipschitzian with the constant $L=K=1$ and $g$ is
Riemann integrable on $\left[ a,b\right] .$

We obviously have 
\begin{align*}
\left| T\left( f,g;u\right) \right| & =\frac{1}{b-a}\int_{a}^{b}f\left(
t\right) g\left( t\right) dt-\frac{1}{b-a}\int_{a}^{b}f\left( t\right)
dt\cdot \frac{1}{b-a}\int_{a}^{b}g\left( t\right) dt \\
& =\frac{b-a}{4},
\end{align*}
\begin{equation*}
\int_{a}^{b}\left| t-\frac{a+b}{2}\right| \left| g\left( t\right) -\frac{1}{%
u\left( b\right) -u\left( a\right) }\int_{a}^{b}g\left( s\right) du\left(
s\right) \right| dt=\frac{\left( b-a\right) ^{2}}{4}
\end{equation*}
\begin{equation*}
\left\| g-\dfrac{1}{u\left( b\right) -u\left( a\right) }\dint_{a}^{b}g\left(
s\right) du\left( s\right) \right\| _{\infty }=\left\| g\right\| _{\infty }=1
\end{equation*}
and 
\begin{equation*}
\left\| g-\dfrac{1}{u\left( b\right) -u\left( a\right) }\dint_{a}^{b}g\left(
s\right) du\left( s\right) \right\| _{p}=\left\| g\right\| _{p}=\left(
b-a\right) ^{\frac{1}{p}}.
\end{equation*}
Consequently, by (\ref{2.14}), one has 
\begin{equation*}
\frac{b-a}{4}\leq \frac{C}{b-a}\frac{\left( b-a\right) ^{2}}{4}\leq \left\{ 
\begin{array}{l}
\dfrac{D\left( b-a\right) ^{2}}{b-a}\cdot 1 \\ 
\\ 
\dfrac{E\left( b-a\right) ^{2}}{\left( q+1\right) ^{\frac{1}{q}}\left(
b-a\right) }%
\end{array}
\right.
\end{equation*}
giving 
\begin{equation*}
\frac{1}{4}\leq \frac{C}{4}\leq \left\{ 
\begin{array}{l}
D \\ 
\\ 
\dfrac{E}{\left( q+1\right) ^{\frac{1}{q}}},\;\;q>1.%
\end{array}
\right.
\end{equation*}
From the first inequality we obtain $C\geq 1.$ Also, we get $D\geq \frac{1}{4%
}$ and $E\geq \frac{\left( q+1\right) ^{\frac{1}{q}}}{4}.$ Letting $%
q\rightarrow 1+,$ we deduce $E\geq \frac{1}{2}$ and the corollary is proved.
\end{proof}

\section{A Quadrature Formula}

Let us consider the partition of the interval $\left[ a,b\right] $ given by 
\begin{equation}
I_{n}:a=x_{0}<x_{1}<\cdots <x_{n-1}<x_{n}=b.  \label{3.1}
\end{equation}
Denote $v\left( I_{n}\right) :=\max \left\{ h_{i}|i=\overline{0,n-1}\right\} 
$ where $h_{i}:=x_{i+1}-x_{i},$ $i=\overline{0,n-1}.$

If $f:\left[ a,b\right] \rightarrow \mathbb{R}$ is continuous on $\left[ a,b%
\right] $ and if we define 
\begin{align*}
M_{i}& :=\sup\limits_{t\in \left[ x_{i},x_{i+1}\right] }f\left( t\right)
,\;\;m_{i}:=\inf\limits_{t\in \left[ x_{i},x_{i+1}\right] }f\left( t\right)
,\;\text{and} \\
v\left( f,I_{n}\right) & =\max\limits_{i=\overline{0,n-1}}\left(
M_{i}-m_{i}\right) ,
\end{align*}
then, obviously, by the continuity of $f$ on $\left[ a,b\right] ,$ for any $%
\varepsilon >0,$ we may find a division $I_{n}$ with norm $v\left(
I_{n}\right) <\delta $ such that $v\left( f,I_{n}\right) <\varepsilon .$

Consider now the quadrature rule 
\begin{equation}
S_{n}\left( f,g;u,I_{n}\right) :=\sum_{i=0}^{n-1}\frac{1}{u\left(
x_{i+1}\right) -u\left( x_{i}\right) }\int_{x_{i}}^{x_{i+1}}f\left( t\right)
du\left( t\right) \cdot \int_{x_{i}}^{x_{i+1}}g\left( t\right) du\left(
t\right)  \label{3.2}
\end{equation}
provided $f,g\in C\left[ a,b\right] ,$ $u\in BV\left[ a,b\right] $ and $%
u\left( x_{i+1}\right) \neq u\left( x_{i}\right) ,$ $i=0,\dots ,n-1.$

We may now state the following result in approximating the Stieltjes
integral 
\begin{equation*}
\int_{a}^{b}f\left( t\right) g\left( t\right) du\left( t\right) .
\end{equation*}

\begin{theorem}
\label{t3.2a}Let $f,g\in C\left[ a,b\right] $ and $u\in BV\left[ a,b\right]
. $ If $I_{n}$ is a division of the interval $\left[ a,b\right] $ and $%
u\left( x_{i+1}\right) \neq u\left( x_{i}\right) ,$ $i=0,\dots ,n-1,$ then
we have: 
\begin{equation}
\int_{a}^{b}f\left( t\right) g\left( t\right) du\left( t\right) =S_{n}\left(
f,g;u,I_{n}\right) +R_{n}\left( f,g;u,I_{n}\right) ,  \label{3.3a}
\end{equation}
where $S_{n}\left( f,g;u,I_{n}\right) $ is as defined in (\ref{3.2}) and the
remainder $R_{n}\left( f,g;u,I_{n}\right) $ satisfies the estimate 
\begin{multline}
\left| R_{n}\left( f,g;u,I_{n}\right) \right| \leq \frac{1}{2}v\left(
f,I_{n}\right)  \label{3.4a} \\
\times \max\limits_{i=\overline{0,n-1}}\left\| g-\frac{1}{u\left(
x_{i+1}\right) -u\left( x_{i}\right) }\int_{x_{i}}^{x_{i+1}}g\left( s\right)
du\left( s\right) \right\| _{\left[ x_{i},x_{i+1}\right] ,\infty
}\bigvee_{a}^{b}\left( u\right) .
\end{multline}
The constant $\frac{1}{2}$ is sharp in (\ref{3.4a}) in the sense that it
cannot be replaced by a smaller constant.
\end{theorem}

\begin{proof}
Applying the inequality (\ref{2.2a}) on the intervals $\left[ x_{i},x_{i+1}%
\right] ,$ $i=0,\dots ,n-1,$ we have 
\begin{multline}
\left| \int_{x_{i}}^{x_{i+1}}f\left( t\right) g\left( t\right) du\left(
t\right) \right.  \label{3.5a} \\
-\left. \frac{1}{u\left( x_{i+1}\right) -u\left( x_{i}\right) }%
\int_{x_{i}}^{x_{i+1}}f\left( t\right) du\left( t\right) \cdot
\int_{x_{i}}^{x_{i+1}}g\left( t\right) du\left( t\right) \right| \\
\leq \frac{1}{2}\left( M_{i}-m_{i}\right) \sup\limits_{t\in \left[
x_{i},x_{i+1}\right] }\left| g\left( t\right) -\frac{1}{u\left(
x_{i+1}\right) -u\left( x_{i}\right) }\int_{x_{i}}^{x_{i+1}}g\left( s\right)
du\left( s\right) \right| \bigvee_{x_{i}}^{x_{i+1}}\left( u\right) .
\end{multline}
Summing the inequalities (\ref{3.5a}) over $i$ from $0$ to $n-1,$ and using
the generalised triangle inequality, we have 
\begin{align}
& \left| R_{n}\left( f,g;u,I_{n}\right) \right|  \label{3.6a} \\
& \leq \frac{1}{2}\sum_{i=0}^{n-1}\left( M_{i}-m_{i}\right) \left\| g-\frac{1%
}{u\left( x_{i+1}\right) -u\left( x_{i}\right) }\int_{x_{i}}^{x_{i+1}}g%
\left( s\right) du\left( s\right) \right\| _{\left[ x_{i},x_{i+1}\right]
,\infty }  \notag \\
&
\;\;\;\;\;\;\;\;\;\;\;\;\;\;\;\;\;\;\;\;\;\;\;\;\;\;\;\;\;\;\;\;\;\;\;\;\;\;%
\;\;\;\;\;\;\;\;\;\;\;\;\;\times \bigvee_{x_{i}}^{x_{i+1}}\left( u\right) 
\notag \\
& \leq \frac{1}{2}v\left( f,I_{n}\right) \max\limits_{i=\overline{0,n-1}%
}\left\| g-\frac{1}{u\left( x_{i+1}\right) -u\left( x_{i}\right) }%
\int_{x_{i}}^{x_{i+1}}g\left( s\right) du\left( s\right) \right\| _{\left[
x_{i},x_{i+1}\right] ,\infty }  \notag \\
&
\;\;\;\;\;\;\;\;\;\;\;\;\;\;\;\;\;\;\;\;\;\;\;\;\;\;\;\;\;\;\;\;\;\;\;\;\;\;%
\;\;\;\;\;\;\;\;\;\;\;\hfill \times
\sum_{i=0}^{n-1}\bigvee_{x_{i}}^{x_{i+1}}\left( u\right)  \notag \\
& =\frac{1}{2}v\left( f,I_{n}\right) \max\limits_{i=\overline{0,n-1}}\left\|
g-\frac{1}{u\left( x_{i+1}\right) -u\left( x_{i}\right) }%
\int_{x_{i}}^{x_{i+1}}g\left( s\right) du\left( s\right) \right\| _{\left[
x_{i},x_{i+1}\right] ,\infty }  \notag \\
&
\;\;\;\;\;\;\;\;\;\;\;\;\;\;\;\;\;\;\;\;\;\;\;\;\;\;\;\;\;\;\;\;\;\;\;\;\;\;%
\;\;\;\;\;\;\;\;\;\;\;\hfill \times \bigvee_{a}^{b}\left( u\right) ,  \notag
\end{align}
and the estimate (\ref{3.4a}) is obtained.
\end{proof}

\begin{remark}
Similar results may be stated for either $u$ monotonic or Lipschitzian. We
omit the details.
\end{remark}

We may now state another result in approximating the Stieltjes integral 
\begin{equation*}
\int_{a}^{b}f\left( t\right) g\left( t\right) du\left( t\right) .
\end{equation*}

\begin{theorem}
\label{t3.2}Let $f,g:\left[ a,b\right] \rightarrow \mathbb{R}$ be such that $%
f$ is of $r-H-$H\"{o}lder type on $\left[ a,b\right] $ (see Theorem \ref%
{t2.1}), $g$ is continuous on $\left[ a,b\right] ,$ $I_{n}$ is as above and $%
u:\left[ a,b\right] \rightarrow \mathbb{R}$ is of bounded variation on $%
\left[ a,b\right] $ with $u\left( x_{i+1}\right) \neq u\left( x_{i}\right) ,$
$i=0,\dots ,n-1.$ Then we have the representation 
\begin{equation}
\int_{a}^{b}f\left( t\right) g\left( t\right) du\left( t\right) =S_{n}\left(
f,g;u,I_{n}\right) +R_{n}\left( f,g;u,I_{n}\right) ,  \label{3.3}
\end{equation}
where the quadrature $S_{n}\left( f,g;u,I_{n}\right) $ is as defined in (\ref%
{3.2}) and the remainder $R_{n}\left( f,g;u,I_{n}\right) $ satisfies the
estimate 
\begin{multline}
\left| R_{n}\left( f,g;u,I_{n}\right) \right| \leq \frac{H}{2^{r}}\left[
v\left( I_{n}\right) \right] ^{r}  \label{3.4} \\
\times \max\limits_{i=\overline{0,n-1}}\left\| g-\frac{1}{u\left(
x_{i+1}\right) -u\left( x_{i}\right) }\int_{x_{i}}^{x_{i+1}}g\left( s\right)
du\left( s\right) \right\| _{\left[ x_{i},x_{i+1}\right] ,\infty
}\bigvee_{a}^{b}\left( u\right) ,
\end{multline}
where $v\left( I_{n}\right) :=\max \left\{ h_{i}|i=\overline{0,n-1}\right\}
. $
\end{theorem}

\begin{proof}
Applying the inequality (\ref{2.2}) on the interval $\left[ x_{i},x_{i+1}%
\right] $ to get 
\begin{multline}
\left| \int_{x_{i}}^{x_{i+1}}f\left( t\right) g\left( t\right) du\left(
t\right) \right.  \label{3.5} \\
-\left. \frac{1}{u\left( x_{i+1}\right) -u\left( x_{i}\right) }%
\int_{x_{i}}^{x_{i+1}}f\left( t\right) du\left( t\right) \cdot
\int_{x_{i}}^{x_{i+1}}g\left( t\right) du\left( t\right) \right| \\
\leq \frac{Hh_{i}^{r}}{2^{r}}\left\| g-\frac{1}{u\left( x_{i+1}\right)
-u\left( x_{i}\right) }\int_{x_{i}}^{x_{i+1}}g\left( t\right) du\left(
t\right) \right\| _{\left[ x_{i},x_{i+1}\right] ,\infty
}\bigvee_{x_{i}}^{x_{i+1}}\left( u\right) ,
\end{multline}
for each $i\in \left\{ 0,\dots ,n-1\right\} .$

Summing the inequalities (\ref{3.5}) over $i$ from $0$ to $n-1,$ and using
the generalised triangle inequality, we have 
\begin{multline}
\left| R_{n}\left( f,g;u,I_{n}\right) \right|  \label{3.6} \\
\leq \frac{H}{2^{r}}\sum_{i=0}^{n-1}h_{i}^{r}\left\| g-\frac{1}{u\left(
x_{i+1}\right) -u\left( x_{i}\right) }\int_{x_{i}}^{x_{i+1}}g\left( t\right)
du\left( t\right) \right\| _{\left[ x_{i},x_{i+1}\right] ,\infty
}\bigvee_{x_{i}}^{x_{i+1}}\left( u\right)
\end{multline}
\begin{align*}
& \leq \frac{H}{2^{r}}\left[ v\left( f\right) \right] ^{n}\max\limits_{i=%
\overline{0,n-1}}\left\| g-\frac{1}{u\left( x_{i+1}\right) -u\left(
x_{i}\right) }\int_{x_{i}}^{x_{i+1}}g\left( t\right) du\left( t\right)
\right\| _{\left[ x_{i},x_{i+1}\right] ,\infty
}\sum_{i=0}^{n-1}\bigvee_{x_{i}}^{x_{i+1}}\left( u\right) \\
& =\frac{H}{2^{r}}\left[ v\left( f\right) \right] ^{n}\max\limits_{i=%
\overline{0,n-1}}\left\| g-\frac{1}{u\left( x_{i+1}\right) -u\left(
x_{i}\right) }\int_{x_{i}}^{x_{i+1}}g\left( s\right) du\left( s\right)
\right\| _{\left[ x_{i},x_{i+1}\right] ,\infty }\bigvee_{a}^{b}\left(
u\right) ,
\end{align*}
and the inequality (\ref{3.4}) is obtained.
\end{proof}

\begin{remark}
Similar results may be stated if one uses Theorem \ref{t2.3} and Theorem \ref%
{t2.5}. We omit the details.
\end{remark}

\section{Some Particular Cases}

For $f,g,w:\left[ a,b\right] \rightarrow \mathbb{R}$, integrable and with
the property that $\int_{a}^{b}w\left( t\right) dt\neq 0,$ reconsider the
weighted \v{C}eby\v{s}ev functional 
\begin{multline}
T_{w}\left( f,g\right) :=\frac{1}{\int_{a}^{b}w\left( t\right) dt}%
\int_{a}^{b}w\left( t\right) f\left( t\right) g\left( t\right) dt
\label{4.1} \\
-\frac{1}{\int_{a}^{b}w\left( t\right) dt}\int_{a}^{b}w\left( t\right)
f\left( t\right) dt\cdot \frac{1}{\int_{a}^{b}w\left( t\right) dt}%
\int_{a}^{b}w\left( t\right) g\left( t\right) dt.
\end{multline}

\noindent \textbf{1. }If $f,g,w:\left[ a,b\right] \rightarrow \mathbb{R}$
are continuous and there exists the real constants $m,M$ such \hspace{0.05in}%
that 
\begin{equation}
m\leq f\left( t\right) \leq M\text{ for each }t\in \left[ a,b\right] ,
\label{4.2}
\end{equation}
then one has the inequality 
\begin{multline}
\left| T_{w}\left( f,g\right) \right| \leq \frac{1}{2}\left( M-m\right) 
\frac{1}{\left| \int_{a}^{b}w\left( s\right) ds\right| }  \label{4.3} \\
\times \left\| g-\frac{1}{\int_{a}^{b}w\left( s\right) ds}%
\int_{a}^{b}g\left( s\right) w\left( s\right) ds\right\| _{\left[ a,b\right]
,\infty }\int_{a}^{b}\left| w\left( s\right) \right| ds.
\end{multline}
The proof follows by Theorem \ref{t2.1a} on choosing $u\left( t\right)
=\int_{a}^{t}w\left( s\right) ds.$

\noindent \textbf{2. }If $f,g,w$ are as in \textbf{1 }and $w\left( s\right)
\geq 0$ for $s\in \left[ a,b\right] ,$ then one has the inequality 
\begin{multline}
\left| T_{w}\left( f,g\right) \right| \leq \frac{1}{2}\left( M-m\right) 
\frac{1}{\int_{a}^{b}w\left( s\right) ds}  \label{4.4a} \\
\times \int_{a}^{b}\left| g\left( t\right) -\frac{1}{\int_{a}^{b}w\left(
s\right) ds}\int_{a}^{b}g\left( s\right) w\left( s\right) ds\right| w\left(
s\right) ds.
\end{multline}
The proof follows by Theorem \ref{t2.2} on choosing $u\left( t\right)
=\int_{a}^{t}w\left( s\right) ds.$

\noindent \textbf{3. }If $f,g$ are Riemann integrable on $\left[ a,b\right] $
and $f$ satisfies (\ref{4.2}), and $w$ is continuous on $\left[ a,b\right] ,$
then one has the inequality 
\begin{multline}
\left| T_{w}\left( f,g\right) \right| \leq \frac{1}{2}\left\| w\right\| _{%
\left[ a,b\right] ,\infty }\left( M-m\right) \frac{1}{\left|
\int_{a}^{b}w\left( s\right) ds\right| }  \label{4.5a} \\
\times \int_{a}^{b}\left| g\left( t\right) -\frac{1}{\int_{a}^{b}w\left(
s\right) ds}\int_{a}^{b}g\left( s\right) w\left( s\right) ds\right| ds.
\end{multline}
The proof follows by Theorem \ref{t2.3} on choosing $u\left( t\right)
=\int_{a}^{t}w\left( s\right) ds.$

\noindent \textbf{4. }If $f,g,w:\left[ a,b\right] \rightarrow \mathbb{R}$
are continuous and $f$ is of $r-H-$H\"{o}lder type (see Theorem \ref{t2.1}),
then one has the inequality 
\begin{multline*}
\left| T_{w}\left( f,g\right) \right| \leq \frac{H\left| b-a\right| ^{r}}{%
2^{r}}\cdot \frac{1}{\left| \int_{a}^{b}w\left( s\right) ds\right| } \\
\times \left\| g-\frac{1}{\int_{a}^{b}w\left( s\right) ds}%
\int_{a}^{b}g\left( s\right) w\left( s\right) ds\right\| _{\left[ a,b\right]
,\infty }\int_{a}^{b}\left| w\left( s\right) \right| ds.
\end{multline*}
The proof follows by Theorem \ref{t2.1} on choosing $u\left( t\right)
=\int_{a}^{t}w\left( s\right) ds.$

\noindent \textbf{5. }If $f,g,w$ are as in \textbf{4 }and $w\left( s\right)
\geq 0$ for $s\in \left[ a,b\right] ,$ then one has the inequality 
\begin{align}
& \left| T_{w}\left( f,g\right) \right|  \label{4.4} \\
& \leq \frac{H}{\int_{a}^{b}w\left( s\right) ds}\int_{a}^{b}\left| t-\frac{%
a+b}{2}\right| ^{r}\left| g\left( t\right) -\frac{1}{\int_{a}^{b}w\left(
s\right) ds}\int_{a}^{b}g\left( s\right) w\left( s\right) ds\right| w\left(
s\right) ds  \notag \\
& \leq \frac{H\left( b-a\right) ^{r}}{2^{r}\int_{a}^{b}w\left( s\right) ds}%
\int_{a}^{b}\left| g\left( t\right) -\frac{1}{\int_{a}^{b}w\left( s\right) ds%
}\int_{a}^{b}g\left( s\right) w\left( s\right) ds\right| w\left( s\right) ds.
\notag
\end{align}
The proof follows by Theorem \ref{t2.3} on choosing $u\left( t\right)
=\int_{a}^{t}w\left( s\right) ds.$

\noindent \textbf{6. }If $f$ is of $r-H-$H\"{o}lder type, $g$ are Riemann
integrable on $\left[ a,b\right] $ and $w$ is continuous on $\left[ a,b%
\right] ,$ then one has the inequality 
\begin{align}
& \left| T_{w}\left( f,g\right) \right|  \label{4.5} \\
& \leq \frac{H\left\| w\right\| _{\left[ a,b\right] ,\infty }}{\left|
\int_{a}^{b}w\left( s\right) ds\right| }\int_{a}^{b}\left| t-\frac{a+b}{2}%
\right| ^{r}\left| g\left( t\right) -\frac{1}{\int_{a}^{b}w\left( s\right) ds%
}\int_{a}^{b}g\left( s\right) w\left( s\right) ds\right| dt  \notag
\end{align}
\begin{equation*}
\leq \left\{ 
\begin{array}{l}
\dfrac{H\left\| w\right\| _{\left[ a,b\right] ,\infty }\left( b-a\right)
^{r+1}}{2^{r}\left( r+1\right) \left| \int_{a}^{b}w\left( s\right) ds\right| 
}\left\| g-\dfrac{1}{\int_{a}^{b}w\left( s\right) ds}\int_{a}^{b}g\left(
s\right) w\left( s\right) ds\right\| _{\left[ a,b\right] ,\infty }; \\ 
\\ 
\dfrac{H\left\| w\right\| _{\left[ a,b\right] ,\infty }\left( b-a\right) ^{r+%
\frac{1}{q}}}{2^{r}\left( qr+1\right) ^{\frac{1}{q}}\left|
\int_{a}^{b}w\left( s\right) ds\right| }\left\| g-\dfrac{1}{%
\int_{a}^{b}w\left( s\right) ds}\int_{a}^{b}g\left( s\right) w\left(
s\right) ds\right\| _{\left[ a,b\right] ,p}, \\ 
\hfill p>1,\;\frac{1}{p}+\frac{1}{q}=1; \\ 
\\ 
\dfrac{H\left\| w\right\| _{\left[ a,b\right] ,\infty }\left( b-a\right) ^{r}%
}{2^{r}\left| \int_{a}^{b}w\left( s\right) ds\right| }\left\| g-\dfrac{1}{%
\int_{a}^{b}w\left( s\right) ds}\int_{a}^{b}g\left( s\right) w\left(
s\right) ds\right\| _{\left[ a,b\right] ,1}.%
\end{array}
\right.
\end{equation*}
The proof follows by Theorem \ref{t2.5} on choosing $u\left( t\right)
=\int_{a}^{t}w\left( s\right) ds.$

\section{Other Inequalities for Stieltjes Integral}

In \cite{1ab}, the authors have considered the following functional 
\begin{equation*}
D\left( f;u\right) :=\int_{a}^{b}f\left( x\right) du\left( x\right) -\left[
u\left( b\right) -u\left( a\right) \right] \cdot \frac{1}{b-a}%
\int_{a}^{b}f\left( t\right) dt,
\end{equation*}
provided that the involved integrals exist.

In the same paper, the following result in estimating the above functional
has been obtained.

\begin{theorem}
\label{ta.1}Let $f,u:\left[ a,b\right] \rightarrow \mathbb{R}$ be such that $%
u$ is Lipschitzian on $\left[ a,b\right] ,$ i.e., 
\begin{equation}
\left| u\left( x\right) -u\left( y\right) \right| \leq L\left| x-y\right| 
\text{ \hspace{0.05in}for any }x,y\in \left[ a,b\right] \;\;\left( L>0\right)
\label{a.1}
\end{equation}
and $f$ is Riemann integrable on $\left[ a,b\right] .$ If $m,M\in \mathbb{R}$
are such that 
\begin{equation}
m\leq f\left( x\right) \leq M\text{ \hspace{0.05in}for any }x,y\in \left[ a,b%
\right] ,  \label{a.2}
\end{equation}
then we have the inequality 
\begin{equation}
\left| D\left( f;u\right) \right| \leq \frac{1}{2}L\left( M-m\right) \left(
b-a\right) .  \label{a.3}
\end{equation}
The constant $\frac{1}{2}$ is sharp in the sense that it cannot be replaced
by a smaller constant.
\end{theorem}

In \cite{1bb}, the following result complementing the above one was obtained.

\begin{theorem}
\label{ta.2}Let $f,u:\left[ a,b\right] \rightarrow \mathbb{R}$ be such that $%
u:\left[ a,b\right] \rightarrow \mathbb{R}$ is of bounded variation in $%
\left[ a,b\right] $ and $f:\left[ a,b\right] \rightarrow \mathbb{R}$ is $K-$%
Lipschitzian $\left( K>0\right) .$ Then we have the inequality 
\begin{equation}
\left\vert D\left( f;u\right) \right\vert \leq \frac{1}{2}K\left( b-a\right)
\bigvee_{a}^{b}\left( u\right) .  \label{a.4}
\end{equation}%
The constant $\frac{1}{2}$ is sharp in the above sense.
\end{theorem}

In this section further similar results will be pointed out. 

The following identity is interesting in itself.

\begin{lemma}
\label{la.3}Let $f,u:\left[ a,b\right] \rightarrow \mathbb{R}$ be such that
the Stieltjes integral $\int_{a}^{b}f\left( t\right) du\left( t\right) $ and
the Riemann integral $\int_{a}^{b}f\left( t\right) dt$ exist. Then we have
the identity 
\begin{align}
D\left( f;u\right) & =\frac{1}{b-a}\int_{a}^{b}\Phi \left( t\right) df\left(
t\right) =\frac{1}{b-a}\int_{a}^{b}\Gamma \left( t\right) df\left( t\right) 
\label{a.5} \\
& =\frac{1}{b-a}\int_{a}^{b}\left( t-a\right) \left( b-t\right) \Delta
\left( t\right) df\left( t\right) ,  \notag
\end{align}%
where 
\begin{eqnarray*}
\Phi \left( t\right)  &:&=\frac{\left( t-a\right) u\left( b\right) +\left(
b-t\right) u\left( a\right) }{b-t}-u\left( t\right) ,\;\;t\in \lbrack a,b),
\\
\Gamma \left( t\right)  &:&=\left( t-a\right) \left[ u\left( b\right)
-u\left( t\right) \right] -\left( b-t\right) \left[ u\left( t\right)
-u\left( a\right) \right] ,\;\;t\in \left[ a,b\right] ,
\end{eqnarray*}%
and 
\begin{equation*}
\Delta \left( t\right) :=\left[ u;b,t\right] -\left[ u;t,a\right] ,\;\;t\in
\left( a,b\right) ,
\end{equation*}%
where $\left[ u;\alpha ,\beta \right] $ is the divided difference, i.e., we
recall it 
\begin{equation*}
\left[ u;\alpha ,\beta \right] :=\frac{u\left( \alpha \right) -u\left( \beta
\right) }{\alpha -\beta }.
\end{equation*}
\end{lemma}

\begin{proof}
We observe that 
\begin{eqnarray*}
\int_{a}^{b}\Phi \left( t\right) df\left( t\right) &=&\int_{a}^{b}\left[ 
\frac{\left( t-a\right) u\left( b\right) +\left( b-t\right) u\left( a\right) 
}{b-t}-u\left( t\right) \right] df\left( t\right) \\
&=&\left. \left[ \frac{\left( t-a\right) u\left( b\right) +\left( b-t\right)
u\left( a\right) }{b-t}-u\left( t\right) \right] f\left( t\right) \right|
_{a}^{b} \\
&&-\int_{a}^{b}f\left( t\right) d\left[ \frac{\left( t-a\right) u\left(
b\right) +\left( b-t\right) u\left( a\right) }{b-t}-u\left( t\right) \right]
\\
&=&\left[ u\left( b\right) -u\left( b\right) \right] -\left[ u\left(
a\right) -u\left( a\right) \right] -\int_{a}^{b}f\left( t\right) \left[ 
\frac{u\left( b\right) -u\left( a\right) }{b-a}dt-du\left( t\right) \right]
\\
&=&\int_{a}^{b}f\left( t\right) du\left( t\right) -\frac{u\left( b\right)
-u\left( a\right) }{b-a}\int_{a}^{b}f\left( t\right) dt
\end{eqnarray*}
and the first identity in (\ref{a.5}) is proved.

The second and third identities are obvious.
\end{proof}

\begin{remark}
\label{ra.4}If $u$ is an integral, i.e., $u\left( t\right)
=\int_{a}^{t}g\left( s\right) ds,$ then from (\ref{a.5}) we deduce Cerone's
result in \cite{0b} 
\begin{equation}
T\left( f,g\right) =\frac{1}{\left( b-a\right) ^{2}}\int_{a}^{b}\Psi \left(
t\right) df\left( t\right) ,  \label{a.6}
\end{equation}
where 
\begin{align*}
\Psi \left( t\right) & =\frac{t-a}{b-t}\int_{a}^{b}g\left( s\right)
ds-\int_{a}^{t}g\left( s\right) ds\;\;\;\left( t\in \lbrack a,b)\right) \\
& =\left( t-a\right) \int_{t}^{b}g\left( s\right) ds-\left( b-t\right)
\int_{a}^{t}g\left( s\right) ds\;\;\;\left( t\in \lbrack a,b]\right) \\
& =\left( t-a\right) \left( b-t\right) \left[ \frac{\int_{t}^{b}g\left(
s\right) ds}{b-t}-\frac{\int_{a}^{t}g\left( s\right) ds}{t-a}\right]
\;\;\;\left( t\in (a,b)\right) .
\end{align*}
If $w:\left[ a,b\right] \rightarrow \mathbb{R}$ is integrable and $%
\int_{a}^{b}w\left( t\right) dt\neq 0,$ then the choice 
\begin{equation}
u\left( t\right) :=\frac{\int_{a}^{t}w\left( s\right) g\left( s\right) ds}{%
\int_{a}^{t}w\left( s\right) ds},\;\;\;t\in \left[ a,b\right] ,  \label{a.7}
\end{equation}
will produce 
\begin{align*}
D\left( f;u\right) & =\frac{\int_{a}^{b}w\left( s\right) f\left( s\right)
g\left( s\right) ds}{\int_{a}^{b}w\left( s\right) ds}-\frac{%
\int_{a}^{b}w\left( s\right) g\left( s\right) ds}{\int_{a}^{b}w\left(
s\right) ds}\cdot \frac{1}{b-a}\int_{a}^{b}f\left( t\right) dt \\
& =:E\left( f,g;w\right) .
\end{align*}
\end{remark}

The following corollary is thus a natural application of the above Lemma \ref%
{la.3}.

\begin{corollary}
\label{ca.5}If $w,f,g$ are Riemann integrable on $\left[ a,b\right] $ and $%
\int_{a}^{b}w\left( t\right) dt\neq 0,$ then 
\begin{align}
E\left( f,g;w\right) & =\int_{a}^{b}\Phi _{w}\left( t\right) df\left(
t\right) =\frac{1}{b-a}\int_{a}^{b}\Gamma _{w}\left( t\right) df\left(
t\right)  \label{a.8} \\
& =\frac{1}{b-a}\int_{a}^{b}\left( t-a\right) \left( b-t\right) \Delta
_{w}\left( t\right) df\left( t\right) ,  \notag
\end{align}
where 
\begin{align*}
\Phi _{w}\left( t\right) & =\left( \frac{t-a}{b-t}\right) \cdot \frac{%
\int_{a}^{b}w\left( s\right) g\left( s\right) ds}{\int_{a}^{b}w\left(
s\right) ds}-\frac{\int_{a}^{t}w\left( s\right) g\left( s\right) ds}{%
\int_{a}^{b}w\left( s\right) ds}, \\
\Gamma _{w}\left( t\right) & =\left( t-a\right) \frac{\int_{a}^{b}w\left(
s\right) g\left( s\right) ds}{\int_{a}^{b}w\left( s\right) ds}-\left(
b-t\right) \frac{\int_{a}^{t}w\left( s\right) g\left( s\right) ds}{%
\int_{a}^{b}w\left( s\right) ds}, \\
\Delta _{w}\left( t\right) & =\frac{\int_{t}^{b}w\left( s\right) g\left(
s\right) ds}{\left( b-t\right) \int_{a}^{b}w\left( s\right) ds}-\frac{%
\int_{a}^{t}w\left( s\right) g\left( s\right) ds}{\left( t-a\right)
\int_{a}^{b}w\left( s\right) ds}.
\end{align*}
\end{corollary}

The following general result in bounding the functional $D\left( f;u\right) $
may be stated.

\begin{theorem}
\label{ta.6}Let $f,u:\left[ a,b\right] \rightarrow \mathbb{R}$.

\begin{enumerate}
\item[$\left( i\right) $] If $f$ is of bounded variation and $u$ is
continuous on $\left[ a,b\right] ,$ then 
\begin{equation}
\left| D\left( f;u\right) \right| \leq \left\{ 
\begin{array}{l}
\sup\limits_{t\in \left[ a,b\right) }\left| \Phi \left( t\right) \right|
\bigvee_{a}^{b}\left( f\right) , \\ 
\\ 
\frac{1}{b-a}\sup\limits_{t\in \left[ a,b\right] }\left| \Gamma \left(
t\right) \right| \bigvee_{a}^{b}\left( f\right) , \\ 
\\ 
\frac{1}{b-a}\sup\limits_{t\in \left( a,b\right) }\left[ \left( t-a\right)
\left( b-t\right) \left| \Delta \left( t\right) \right| \right]
\bigvee_{a}^{b}\left( f\right) .%
\end{array}
\right.  \label{a.9}
\end{equation}

\item[$\left( ii\right) $] If $f$ is $L-$Lipschitzian and $u$ is Riemann
integrable on $\left[ a,b\right] ,$ then 
\begin{equation}
\left| D\left( f;u\right) \right| \leq \left\{ 
\begin{array}{l}
L\int_{a}^{b}\left| \Phi \left( t\right) \right| dt, \\ 
\\ 
\frac{L}{b-a}\int_{a}^{b}\left| \Gamma \left( t\right) \right| dt, \\ 
\\ 
\frac{L}{b-a}\int_{a}^{b}\left( t-a\right) \left( b-t\right) \left| \Delta
\left( t\right) \right| dt.%
\end{array}
\right.  \label{a.10}
\end{equation}

\item[$\left( iii\right) $] If $f$ is monotonic nondecreasing on $\left[ a,b%
\right] $ and $u$ is continuous on $\left[ a,b\right] ,$ then 
\begin{equation}
\left| D\left( f;u\right) \right| \leq \left\{ 
\begin{array}{l}
\int_{a}^{b}\left| \Phi \left( t\right) \right| df\left( t\right) , \\ 
\\ 
\frac{1}{b-a}\int_{a}^{b}\left| \Gamma \left( t\right) \right| df\left(
t\right) , \\ 
\\ 
\frac{1}{b-a}\int_{a}^{b}\left( t-a\right) \left( b-t\right) \left| \Delta
\left( t\right) \right| df\left( t\right) .%
\end{array}
\right.  \label{a.11}
\end{equation}
\end{enumerate}
\end{theorem}

\begin{proof}
Follows by Lemma \ref{la.3} on taking into account that 
\begin{equation*}
\left| \int_{c}^{d}p\left( t\right) dv\left( t\right) \right| \leq
\sup\limits_{t\in \left[ a,b\right] }\left| p\left( t\right) \right|
\bigvee_{c}^{d}\left( t\right)
\end{equation*}
if $p$ is continuous on $\left[ c,d\right] $ and $v$ is of bounded
variation, 
\begin{equation*}
\left| \int_{c}^{d}p\left( t\right) dv\left( t\right) \right| \leq
L\int_{c}^{d}\left| p\left( t\right) \right| dt;
\end{equation*}
if $v$ is $L-$Lipschitzian on $\left[ c,d\right] $ and $p$ is Riemann
integrable on $\left[ c,d\right] $ and 
\begin{equation*}
\left| \int_{c}^{d}p\left( t\right) dv\left( t\right) \right| \leq
\int_{c}^{d}\left| p\left( t\right) \right| dt
\end{equation*}
if $p$ is continuous on $\left[ c,d\right] $ and $v$ is monotonic
nondecreasing on $\left[ c,d\right] .$
\end{proof}

It is natural to consider the following corollaries, since they provide
simpler bounds for the functional $D\left( f;u\right) $ in terms of $\Delta $
defined in Lemma \ref{la.3}.

\begin{corollary}
\label{ca.7}If $f$ is of bounded variation and $u$ is continuous on $\left[
a,b\right] ,$ then 
\begin{align}
\left| D\left( f;u\right) \right| & \leq \frac{1}{b-a}\sup\limits_{t\in %
\left[ a,b\right] }\left[ \left( t-a\right) \left( b-t\right) \Delta \left(
t\right) \right] \bigvee_{a}^{b}\left( f\right)  \label{a.12} \\
& \leq \frac{b-a}{4}\left\| \Delta \right\| _{\infty }\bigvee_{a}^{b}\left(
f\right) .  \notag
\end{align}
\end{corollary}

\begin{corollary}
\label{ca.8}If $f$ is $L-$Lipschitzian and $u$ is Riemann integrable on $%
\left[ a,b\right] ,$ then 
\begin{eqnarray}
&&\left\vert D\left( f;u\right) \right\vert   \label{a.13} \\
&\leq &\frac{L}{b-a}\int_{a}^{b}\left( t-a\right) \left( b-t\right)
\left\vert \Delta \left( t\right) \right\vert dt \\
&\leq &\left\{ 
\begin{array}{l}
\frac{1}{6}L\left( b-a\right) ^{2}\left\Vert \Delta \right\Vert _{\infty },
\\ 
\\ 
L\left( b-a\right) ^{1+\frac{1}{q}}\left[ B\left( q+1,q+1\right) \right] ^{%
\frac{1}{q}}\left\Vert \Delta \right\Vert _{p},\;p>1,\;\frac{1}{p}+\frac{1}{q%
}=1; \\ 
\\ 
\frac{1}{4}L\left( b-a\right) \left\Vert \Delta \right\Vert _{1},%
\end{array}%
\right.   \notag
\end{eqnarray}%
and 
\begin{equation*}
\end{equation*}
\end{corollary}

\begin{corollary}
\label{ca.9}If $f$ is monotonic nondecreasing and $g$ is continuous, then 
\begin{align}
& \left\vert D\left( f;u\right) \right\vert   \label{a.14} \\
& \leq \frac{1}{b-a}\int_{a}^{b}\left( t-a\right) \left( b-t\right)
\left\vert \Delta \left( t\right) \right\vert dt  \notag \\
& \leq \left\{ 
\begin{array}{l}
\frac{1}{4}\left( b-a\right) \int_{a}^{b}\left\vert \Delta \left( t\right)
\right\vert df\left( t\right) , \\ 
\\ 
\frac{1}{b-a}\left( \int_{a}^{b}\left[ \left( b-t\right) \left( t-a\right) %
\right] ^{q}df\left( t\right) \right) ^{\frac{1}{q}}\left(
\int_{a}^{b}\left\vert \Delta \left( t\right) \right\vert ^{p}df\left(
t\right) \right) ^{\frac{1}{p}}, \\ 
\hfill p>1,\;\frac{1}{p}+\frac{1}{q}=1; \\ 
\frac{1}{b-a}\left\Vert \Delta \right\Vert _{\infty }\int_{a}^{b}\left(
t-a\right) \left( b-t\right) df\left( t\right) .%
\end{array}%
\right.   \notag
\end{align}
\end{corollary}

\begin{remark}
\label{ra.10}If one chooses in Corollaries \ref{ca.7} -- \ref{ca.9}, $%
u\left( t\right) =\int_{a}^{t}g\left( s\right) ds,$ then the result
incorporated in Theorems 4 -- 6 of \cite{0b} are recaptured.
\end{remark}

Finally, the following result on the positivity of the functional $D\left(
f;u\right) $ holds.

\begin{theorem}
\label{ta.11}Let $f$ be a monotonic nondecreasing function on $\left[ a,b%
\right] .$ If $u$ is such that 
\begin{equation}
\Delta \left( t\right) =\Delta \left( u;a,t,b\right) :=\left[ u;b,t\right] -%
\left[ u;t,a\right] \geq 0  \label{a.15}
\end{equation}%
for any $t\in \left( a,b\right) ,$ then we have the inequality 
\begin{equation}
D\left( f;u\right) \geq \frac{1}{b-a}\left\vert \int_{a}^{b}\left(
t-a\right) \left( b-t\right) \left[ \left\vert \left[ u;b,t\right]
\right\vert -\left\vert \left[ u;t,a\right] \right\vert \right] df\left(
t\right) \right\vert \geq 0.  \label{a.16}
\end{equation}
\end{theorem}

The proof is similar to the case in Theorem 3 of \cite{CD} and we omit the
details.

\begin{remark}
\label{ra.11}It is easy to see that a sufficient condition for (\ref{a.15})
to hold is that $u:\left[ a,b\right] \rightarrow \mathbb{R}$ is a convex
function on $\left[ a,b\right] .$
\end{remark}

\begin{remark}
\label{ra.12}Similar results for composite rules in approximating the
Stieltjes integral may be stated but we omit the details.
\end{remark}

\section{Some Results for Monotonic Integrators}

The following result holds.

\begin{theorem}
\label{tb.1}Let $f:\left[ a,b\right] \rightarrow \mathbb{R}$ be $L-$%
Lipschitzian on $\left[ a,b\right] $ and $u$ monotonic nondecreasing on $%
\left[ a,b\right] .$ Then we have he inequality 
\begin{align}
\left\vert D\left( f;u\right) \right\vert & \leq \frac{1}{2}L\left(
b-a\right) \left[ u\left( b\right) -u\left( a\right) -K\left( u\right) %
\right]   \label{b.1} \\
& \leq \frac{1}{2}L\left( b-a\right) \left[ u\left( b\right) -u\left(
a\right) \right] ,  \notag
\end{align}%
where 
\begin{align}
K\left( u\right) & :=\frac{4}{\left( b-a\right) ^{2}}\int_{a}^{b}u\left(
x\right) \left( x-\frac{a+b}{2}\right) dx  \label{b.2} \\
& =\frac{4}{\left( b-a\right) ^{2}}\int_{a}^{b}\left[ u\left( x\right)
-u\left( \frac{a+b}{2}\right) \right] \left( x-\frac{a+b}{2}\right) dx\geq 0.
\notag
\end{align}%
The constant $\frac{1}{2}$ in both inequalities is sharp in the sense that
it cannot be replaced by a smaller constat.
\end{theorem}

\begin{proof}
As $u$ is monotonic nondecreasing on $\left[ a,b\right] ,$ then 
\begin{align}
& \left\vert \int_{a}^{b}f\left( x\right) du\left( x\right) -\frac{u\left(
b\right) -u\left( a\right) }{b-a}\int_{a}^{b}f\left( t\right) dt\right\vert 
\label{b.3} \\
& =\left\vert \int_{a}^{b}\left( f\left( x\right) -\frac{1}{b-a}%
\int_{a}^{b}f\left( t\right) dt\right) du\left( x\right) \right\vert   \notag
\\
& \leq \int_{a}^{b}\left\vert f\left( x\right) -\frac{1}{b-a}%
\int_{a}^{b}f\left( t\right) dt\right\vert du\left( x\right) .  \notag
\end{align}%
Taking into account that $f$ is $L-$Lipschitzian, we have the following
Ostrowski type inequality (see for example \cite{SEV1}) 
\begin{equation}
\left\vert f\left( x\right) -\frac{1}{b-a}\int_{a}^{b}f\left( t\right)
dt\right\vert \leq L\left[ \frac{1}{4}+\left( \frac{x-\frac{a+b}{2}}{b-a}%
\right) ^{2}\right] \left( b-a\right)   \label{b.4}
\end{equation}%
for all $x\in \left[ a,b\right] ,$ from where we deduce 
\begin{multline}
\int_{a}^{b}\left\vert f\left( x\right) -\frac{1}{b-a}\int_{a}^{b}f\left(
t\right) dt\right\vert du\left( x\right)   \label{b.5} \\
\leq L\left( b-a\right) \int_{a}^{b}\left[ \frac{1}{4}+\left( \frac{x-\frac{%
a+b}{2}}{b-a}\right) ^{2}\right] du\left( x\right) .
\end{multline}%
Now, observe that, by the integration by parts formula for the Stieltjes
integral, we have 
\begin{eqnarray*}
\int_{a}^{b}\left( x-\frac{a+b}{2}\right) ^{2}du\left( x\right)  &=&u\left(
x\right) \left. \left( x-\frac{a+b}{2}\right) ^{2}\right\vert
_{a}^{b}-2\int_{a}^{b}u\left( x\right) \left( x-\frac{a+b}{2}\right) dx \\
&=&\frac{\left( b-a\right) ^{2}}{4}\left[ u\left( b\right) -u\left( a\right) %
\right] -2\int_{a}^{b}u\left( x\right) \left( x-\frac{a+b}{2}\right) dx
\end{eqnarray*}%
and then 
\begin{multline}
\int_{a}^{b}\left[ \frac{1}{4}+\left( \frac{x-\frac{a+b}{2}}{b-a}\right) ^{2}%
\right] du\left( x\right)   \label{b.6} \\
=\frac{1}{2}\left[ u\left( b\right) -u\left( a\right) \right] -\frac{2}{%
\left( b-a\right) ^{2}}\int_{a}^{b}u\left( x\right) \left( x-\frac{a+b}{2}%
\right) dx.
\end{multline}%
Using (\ref{b.3}) -- (\ref{b.6}) we deduce the first part of (\ref{b.1}).

The second part is obvious by (\ref{b.2}) which follows by the monotonicity
of $u$ on $\left[ a,b\right] .$

To prove the sharpness of the constant $\frac{1}{2},$ assume that (\ref{b.7}%
) holds with the constants $C,D>0,$ i.e., 
\begin{equation}
\left| D\left( f;u\right) \right| \leq CL\left( b-a\right) \left[ u\left(
b\right) -u\left( a\right) -K\left( u\right) \right] \leq DL\left(
b-a\right) \left[ u\left( b\right) -u\left( a\right) \right] .  \label{b.7}
\end{equation}
Consider the functions $f,u:\left[ a,b\right] \rightarrow \mathbb{R}$ given
by $f\left( x\right) =x-\frac{a+b}{2}$ and 
\begin{equation*}
u\left( x\right) =\left\{ 
\begin{array}{ll}
0 & \text{if \hspace{0.05in}}x\in \lbrack a,b) \\ 
&  \\ 
1 & \text{if \hspace{0.05in}}x=b.%
\end{array}
\right.
\end{equation*}
Thus $f$ is $L-$Lipschitzian with the constant $L=1$ and $u$ is monotonic
nondecreasing.

We observe that 
\begin{equation*}
D\left( f;u\right) =\int_{a}^{b}f\left( x\right) du\left( x\right) =f\left(
x\right) u\left( x\right) \bigg|_{a}^{b}-\int_{a}^{b}u\left( x\right) dx=%
\frac{b-a}{2},
\end{equation*}
$K\left( u\right) =0$ and $u\left( b\right) -u\left( a\right) =1,$ giving in
(\ref{b.7}) 
\begin{equation*}
\frac{b-a}{2}\leq C\left( b-a\right) \leq D\left( b-a\right)
\end{equation*}
and thus $C,D\geq \frac{1}{2}$ proving the sharpness of the constant $\frac{1%
}{2}$ in (\ref{b.1}).
\end{proof}

Another result of this type is the following one.

\begin{theorem}
\label{tb.2}Let $u:\left[ a,b\right] \rightarrow \mathbb{R}$ be monotonic
nondecreasing on $\left[ a,b\right] $ and $f:\left[ a,b\right] \rightarrow 
\mathbb{R}$ be of bounded variation such that the Stieltjes integral $%
\int_{a}^{b}f\left( x\right) du\left( x\right) $ exists. Then we have the
inequality 
\begin{align}
\left\vert D\left( f;u\right) \right\vert & \leq \left[ u\left( b\right)
-u\left( a\right) -Q\left( u\right) \right] \bigvee_{a}^{b}\left( f\right) 
\label{b.8} \\
& \leq \left[ u\left( b\right) -u\left( a\right) \right] \bigvee_{a}^{b}%
\left( f\right) ,  \notag
\end{align}%
where 
\begin{align}
Q\left( u\right) & :=\frac{1}{b-a}\int_{a}^{b}\func{sgn}\left( x-\frac{a+b}{2%
}\right) u\left( x\right) dx  \label{b.8.a} \\
& =\frac{1}{b-a}\int_{a}^{b}\func{sgn}\left( x-\frac{a+b}{2}\right) \left[
u\left( x\right) -u\left( \frac{a+b}{2}\right) \right] dx\geq 0.  \notag
\end{align}%
The first inequality in (\ref{b.8}) is sharp in the sense that the constant $%
c=1$ cannot be replaced by a smaller constant.
\end{theorem}

\begin{proof}
Since $u$ is monotonic nondecreasing, we have (see (\ref{b.3})) that 
\begin{equation}
\left\vert D\left( f;u\right) \right\vert \leq \int_{a}^{b}\left\vert
f\left( x\right) -\frac{1}{b-a}\int_{a}^{b}f\left( t\right) dt\right\vert
du\left( x\right) .  \label{b.9}
\end{equation}%
Using the following Ostrowski type inequality obtained by the author in \cite%
{SEV2} 
\begin{equation}
\left\vert f\left( x\right) -\frac{1}{b-a}\int_{a}^{b}f\left( t\right)
dt\right\vert \leq \left[ \frac{1}{2}+\left\vert \frac{x-\frac{a+b}{2}}{b-a}%
\right\vert \right] \bigvee_{a}^{b}\left( f\right)   \label{b.10}
\end{equation}%
for any $x\in \left[ a,b\right] ,$ we have 
\begin{equation}
\int_{a}^{b}\left\vert f\left( x\right) -\frac{1}{b-a}\int_{a}^{b}f\left(
t\right) dt\right\vert du\left( x\right) \leq \bigvee_{a}^{b}\left( f\right)
\int_{a}^{b}\left[ \frac{1}{2}+\left\vert \frac{x-\frac{a+b}{2}}{b-a}%
\right\vert \right] du\left( x\right) .  \label{b.11}
\end{equation}%
A simple calculation with the Stieltjes integral gives that 
\begin{eqnarray}
&&\int_{a}^{b}\left\vert x-\frac{a+b}{2}\right\vert du\left( x\right) 
\label{b.12} \\
&=&\int_{a}^{\frac{a+b}{2}}\left( \frac{a+b}{2}-x\right) du\left( x\right)
+\int_{\frac{a+b}{2}}^{b}\left( x-\frac{a+b}{2}\right) du\left( x\right)  
\notag \\
&=&\left. u\left( x\right) \left( \frac{a+b}{2}-x\right) \right\vert _{a}^{%
\frac{a+b}{2}}+\int_{a}^{\frac{a+b}{2}}u\left( x\right) dx  \notag \\
&&+\left. \left( x-\frac{a+b}{2}\right) u\left( x\right) \right\vert _{\frac{%
a+b}{2}}^{b}-\int_{\frac{a+b}{2}}^{b}u\left( x\right) dx  \notag \\
&=&\frac{1}{2}\left( b-a\right) \left[ u\left( b\right) -u\left( a\right) %
\right] -\int_{a}^{b}\func{sgn}\left( x-\frac{a+b}{2}\right) u\left(
x\right) dx  \notag
\end{eqnarray}%
and then by (\ref{b.9}) -- (\ref{b.12}) we deduce the first inequality in (%
\ref{b.8}).

The second part of (\ref{b.8}) follows by (\ref{b.8.a}) which holds by the
monotonicity property of $u.$

Now, assume that the first inequality in (\ref{b.8}) holds with a constant $%
E>0,$ i.e., 
\begin{equation}
\left| D\left( f;u\right) \right| \leq E\bigvee_{a}^{b}\left( f\right) \left[
u\left( b\right) -u\left( a\right) -Q\left( u\right) \right] .  \label{b.13}
\end{equation}
Consider the mappings $f,u:\left[ a,b\right] \rightarrow \mathbb{R}$, $%
f\left( x\right) =x-\frac{a+b}{2},$ and 
\begin{equation*}
u\left( x\right) =\left\{ 
\begin{array}{ll}
0 & \text{if \hspace{0.05in}}x\in \left[ a,\frac{a+b}{2}\right] , \\ 
&  \\ 
1 & \text{if \hspace{0.05in}}x\in \left( \frac{a+b}{2},b\right] .%
\end{array}
\right.
\end{equation*}
Then we have 
\begin{eqnarray*}
D\left( f;u\right) &=&\int_{a}^{b}f\left( x\right) du\left( x\right) -\frac{%
u\left( b\right) -u\left( a\right) }{b-a}\int_{a}^{b}f\left( t\right) dt \\
&=&\int_{a}^{b}\left( x-\frac{a+b}{2}\right) du\left( x\right) =\left.
\left( x-\frac{a+b}{2}\right) u\left( x\right) \right|
_{a}^{b}-\int_{a}^{b}u\left( x\right) dx \\
&=&\frac{b-a}{2}\left[ u\left( b\right) +u\left( a\right) \right] =\frac{b-a%
}{2}
\end{eqnarray*}
and 
\begin{eqnarray*}
&&\bigvee_{a}^{b}\left( f\right) \left[ u\left( b\right) -u\left( a\right)
-Q\left( u\right) \right] \\
&=&\left( b-a\right) \left[ u\left( b\right) -u\left( a\right) -\left( \frac{%
1}{b-a}\int_{a}^{\frac{a+b}{2}}\func{sgn}\left( x-\frac{a+b}{2}\right)
u\left( x\right) dx\right. \right. \\
&&+\left. \left. \frac{1}{b-a}\int_{\frac{a+b}{2}}^{b}\func{sgn}\left( x-%
\frac{a+b}{2}\right) u\left( x\right) dx\right) \right] \\
&=&\frac{b-a}{2}.
\end{eqnarray*}
Thus, by (\ref{b.13}) we obtain 
\begin{equation*}
\frac{b-a}{2}\leq E\cdot \frac{b-a}{2},
\end{equation*}
showing that $E\geq 1,$ and the theorem is proved.
\end{proof}

\begin{remark}
\label{rnew}Similar results for composite rules in approximating the
Stieltjes integral may be stated, but we omit the details.
\end{remark}

For other inequalities of  Gr\"{u}ss type, see \cite{Pa1}-\cite{Lup}.

\end{document}